\documentclass[12pt]{amsart}
\usepackage{amsmath}
\usepackage{amscd}
\usepackage{amssymb,latexsym}
\usepackage{epsf}

\newtheorem{theorem}{Theorem}[section]
\newtheorem{proposition}[theorem]{Proposition}

\newtheorem{ex}[theorem]{Example}
\newenvironment{example}{\begin{ex}\rm}{\end{ex}}

\numberwithin{theorem}{section}
\def\Mtwo{{\sl Macaulay 2}}

{\obeyspaces\global\let =\ }
\font\tteight=cmtt10

\copyrightinfo{2000}{(Frank Sottile)}

\begin{document}

\title[Enumerative geometry via Macaulay 2]{An excursion from
enumerative geometry to\\solving systems of polynomial equations\\
with Macaulay 2} 
\author{Frank Sottile}
\address{Department of Mathematics and Statistics\\
        University of Massachusetts\\
        Amherst, MA, 01003\\
        USA}

\email{sottile@math.umass.edu}
\urladdr{http://www.math.umass.edu/\~{}sottile}
\date{14 November 2000}
\thanks{Research supported in part by NSF grant DMS-0070494}
\thanks{Chapter in book: {\it Computations in Algebraic Geometry with
         \Mtwo{}}, D.~Eisenbud, D.~Grayson, M.~Stillmann, and B.~Sturmfels, eds.}

\begin{abstract}
Solving a system of polynomial equations is a ubiquitous problem in
the applications of mathematics. 
Until recently, it has been hopeless to find explicit solutions to such
systems, and mathematics has instead developed deep and
powerful theories about the solutions to polynomial equations.
Enumerative Geometry is concerned with counting the
number of solutions when the polynomials come from a geometric situation and
Intersection Theory gives methods to accomplish the enumeration.

   We use \Mtwo{}\/ to investigate some problems from enumerative geometry,
illustrating some applications of symbolic computation to this important
problem of solving systems of polynomial equations.
Besides enumerating solutions to the resulting polynomial systems, which
include overdetermined, deficient, and improper systems, we address the
important question of real solutions to these geometric problems. 

  The text contains evaluated \Mtwo{}\/ code to illuminate the discussion.
This is intended as a chapter in a book on applications of Macaulay 2 to
problems in mathematics.
While this chapter is largely expository, the results in the last section
concerning lines tangent to quadrics are new.
\end{abstract}

\maketitle

\section{Introduction}
A basic question to ask about a system of polynomial equations is
its number of solutions\index{polynomial equations}.
For this, the fundamental result is the following  
B\'ezout Theorem\index{B\'ezout Theorem}.

\begin{theorem}
 The number of isolated solutions to a system of polynomial equations
$$
 f_1(x_1,\ldots,x_n)=f_2(x_1,\ldots,x_n)= \cdots
 =f_n(x_1,\ldots,x_n)=0
$$
 is bounded by $d_1d_2\cdots d_n$, where $d_i:=\deg f_i$.
 If the polynomials are generic, then this bound is attained for 
 solutions in an algebraically closed field\index{field!algebraically closed}.
\end{theorem}

Here, isolated is taken with respect to the algebraic closure.
This B\'ezout Theorem is a consequence of the refined B\'ezout Theorem of
Fulton and MacPherson~\cite[\S 1.23]{SO:Fu84a}.

A system of polynomial equations with fewer than this
degree bound or B\'ezout number\index{B\'ezout number} of solutions is called
{\it deficient}\index{polynomial equations!deficient}, 
and there are well-defined classes of deficient systems that satisfy other
bounds.  
For example, fewer monomials lead to fewer solutions, for which polyhedral
bounds~\cite{SO:Bernstein} on the number of solutions are often tighter (and
no weaker than) the B\'ezout number, which applies when
all monomials are present.
When the polynomials come from geometry, determining the
number of solutions is the central problem in enumerative
geometry\index{enumerative geometry}.

Symbolic computation\index{symbolic computation} can help compute the
solutions to a system of equations 
that has only isolated solutions.
In this case, the polynomials generate a zero-dimensional ideal
$I$\index{artinian, {\it see also} ideal,
zero-dimensional}\index{ideal!zero-dimensional}. 
The {\it degree}\index{ideal!degree} of $I$ ($\dim_k k[X]/I$), which is the
number of standard monomials in any term order gives an upper bound on
the number of solutions, which is attained when $I$ is
radical\index{ideal!radical}. 

\begin{example}\label{ex:one}
We illustrate this discussion with an example.
Let $f_1$, $f_2$, $f_3$, and $f_4$ be random quadratic polynomials in the ring  
${\mathbb F}_{101}[y_{11},y_{12},y_{21},y_{22}]$.
\par
\vskip 5 pt
\begingroup
\tteight
\baselineskip=10.01pt
\lineskip=0pt
\obeyspaces
i1\ :\ R\ =\ ZZ/101[y11,\ y12,\ y21,\ y22];\leavevmode\hss\endgraf
\endgroup
\penalty-1000
\par
\vskip 1 pt
\noindent
\par
\vskip 5 pt
\begingroup
\tteight
\baselineskip=10.01pt
\lineskip=0pt
\obeyspaces
i2\ :\ PolynomialSystem\ =\ apply(1..4,\ i\ ->\ \leavevmode\hss\endgraf
\ \ \ \ \ \ \ \ \ \ \ \ \ \ \ \ \ \ \ \ \ random(0,\ R)\ +\ random(1,\ R)\ +\ random(2,\ R));\leavevmode\hss\endgraf
\endgroup
\penalty-1000
\par
\vskip 1 pt
\noindent
The ideal they generate has dimension 0 and degree $16=2^4$, which is the
B\'ezout number. 
\par
\vskip 5 pt
\begingroup
\tteight
\baselineskip=10.01pt
\lineskip=0pt
\obeyspaces
i3\ :\ I\ =\ ideal\ PolynomialSystem;\leavevmode\hss\endgraf
\penalty-500\leavevmode\hss\endgraf
o3\ :\ Ideal\ of\ R\leavevmode\hss\endgraf
\endgroup
\penalty-1000
\par
\vskip 1 pt
\noindent
\par
\vskip 5 pt
\begingroup
\tteight
\baselineskip=10.01pt
\lineskip=0pt
\obeyspaces
i4\ :\ dim\ I,\ degree\ I\leavevmode\hss\endgraf
\penalty-500\leavevmode\hss\endgraf
o4\ =\ (0,\ 16)\leavevmode\hss\endgraf
\penalty-500\leavevmode\hss\endgraf
o4\ :\ Sequence\leavevmode\hss\endgraf
\endgroup
\penalty-1000
\par
\vskip 1 pt
\noindent
If we restrict the monomials which appear in the $f_i$ to be among
$$
  1,\;\ y_{11},\;\ y_{12},\;\ y_{21},\;\ y_{22},\;\  
  y_{11}y_{22},\;\ \mbox{ and }\;\  y_{12}y_{21},
$$
then the ideal they generate again has dimension 0, but its degree is now 4.
\par
\vskip 5 pt
\begingroup
\tteight
\baselineskip=10.01pt
\lineskip=0pt
\obeyspaces
i5\ :\ J\ =\ ideal\ (random(R{\char`\^}4,\ R{\char`\^}7)\ *\ \ transpose(\leavevmode\hss\endgraf
\ \ \ \ \ \ \ \ \ \ \ \ \ matrix{\char`\{}{\char`\{}1,\ y11,\ y12,\ y21,\ y22,\ y11*y22,\ y12*y21{\char`\}}{\char`\}}));\leavevmode\hss\endgraf
\penalty-500\leavevmode\hss\endgraf
o5\ :\ Ideal\ of\ R\leavevmode\hss\endgraf
\endgroup
\penalty-1000
\par
\vskip 1 pt
\noindent
\par
\vskip 5 pt
\begingroup
\tteight
\baselineskip=10.01pt
\lineskip=0pt
\obeyspaces
i6\ :\ dim\ J,\ degree\ J\leavevmode\hss\endgraf
\penalty-500\leavevmode\hss\endgraf
o6\ =\ (0,\ 4)\leavevmode\hss\endgraf
\penalty-500\leavevmode\hss\endgraf
o6\ :\ Sequence\leavevmode\hss\endgraf
\endgroup
\penalty-1000
\par
\vskip 1 pt
\noindent
If we further require that the coefficients of the
quadratic terms sum to zero, then the ideal they generate now 
has degree 2.
\par
\vskip 5 pt
\begingroup
\tteight
\baselineskip=10.01pt
\lineskip=0pt
\obeyspaces
i7\ :\ K\ =\ ideal\ (random(R{\char`\^}4,\ R{\char`\^}6)\ *\ transpose(\ \leavevmode\hss\endgraf
\ \ \ \ \ \ \ \ \ \ \ \ \ matrix{\char`\{}{\char`\{}1,\ y11,\ y12,\ y21,\ y22,\ y11*y22\ -\ y12*y21{\char`\}}{\char`\}}));\leavevmode\hss\endgraf
\penalty-500\leavevmode\hss\endgraf
o7\ :\ Ideal\ of\ R\leavevmode\hss\endgraf
\endgroup
\penalty-1000
\par
\vskip 1 pt
\noindent
\par
\vskip 5 pt
\begingroup
\tteight
\baselineskip=10.01pt
\lineskip=0pt
\obeyspaces
i8\ :\ dim\ K,\ degree\ K\leavevmode\hss\endgraf
\penalty-500\leavevmode\hss\endgraf
o8\ =\ (0,\ 2)\leavevmode\hss\endgraf
\penalty-500\leavevmode\hss\endgraf
o8\ :\ Sequence\leavevmode\hss\endgraf
\endgroup
\penalty-1000
\par
\vskip 1 pt
\noindent
In Example~\ref{ex:G22}, we shall see how this last specialization is
geometrically meaningful.
\end{example}

For us, enumerative geometry\index{enumerative geometry} is concerned
with {\sl enumerating geometric figures of some kind having
specified positions with respect to general fixed figures}.
That is, counting the solutions to a geometrically meaningful
system of polynomial equations\index{polynomial equations}.
We use \Mtwo{}\index{\Mtwo{}}\/ to investigate some enumerative geometric
problems\index{enumerative problem} from this point of view.
The problem of enumeration will be solved by computing the
degree\index{ideal!degree} of the 
(0-dimensional) ideal generated by the polynomials.

\section{Solving systems of polynomials}
We briefly discuss some aspects of solving systems of polynomial
equations\index{solving polynomial equations}.
For a more complete survey, see the relevant chapters
in~\cite{SO:CCS,SO:CLO92}.

Given an ideal $I$ in a polynomial ring $k[X]$, set 
${\mathcal V}(I):= {\rm Spec}\,k[X]/I$.
When $I$ is generated by the polynomials
$f_1,\ldots,f_N$, ${\mathcal V}(I)$ gives the set of solutions in affine 
space to the system 
\begin{equation}\label{eq:system}
  f_1(X)\ =\ \cdots\ =\ f_N(X)\ =\ 0
\end{equation}
a geometric structure.
These solutions are the {\it roots} of the ideal $I$.
The degree of a zero-dimensional ideal $I$ provides an algebraic count of
its roots.
The degree of its radical counts 
roots in the algebraic closure, ignoring multiplicities.
          
\subsection{Excess intersection}
Sometimes, only a proper (open) subset of affine space is 
geometrically meaningful, and we want to count only the meaningful roots of
$I$. 
Often the roots ${\mathcal V}(I)$ has positive dimensional components that lie
in the complement of the meaningful subset.
One way to treat this situation of excess or improper intersection is to
saturate\index{saturate} $I$ by a polynomial $f$ vanishing on the extraneous
roots. 
This has the effect of working in $k[X][f^{-1}]$, the coordinate ring of the
complement of ${\mathcal V}(f)$~\cite[Exer.~2.3]{SO:MR97a:13001}.

\begin{example}\label{ex:two}
We illustrate this with an example.
Consider the following ideal in ${\mathbb F}_7[x,y]$.
\par
\vskip 5 pt
\begingroup
\tteight
\baselineskip=10.01pt
\lineskip=0pt
\obeyspaces
i9\ :\ R\ =\ ZZ/7[y,\ x,\ MonomialOrder=>Lex];\leavevmode\hss\endgraf
\endgroup
\penalty-1000
\par
\vskip 1 pt
\noindent
\par
\vskip 5 pt
\begingroup
\tteight
\baselineskip=10.01pt
\lineskip=0pt
\obeyspaces
i10\ :\ I\ =\ ideal\ (y{\char`\^}3*x{\char`\^}2\ +\ 2*y{\char`\^}2*x\ +\ 3*x*y,\ \ 3*y{\char`\^}2\ +\ x*y\ -\ 3*y);\leavevmode\hss\endgraf
\penalty-500\leavevmode\hss\endgraf
o10\ :\ Ideal\ of\ R\leavevmode\hss\endgraf
\endgroup
\penalty-1000
\par
\vskip 1 pt
\noindent
Since the generators have greatest common factor $y$, $I$ defines
finitely many points together with the line $y=0$.
Saturate $I$\/ by the variable $y$ to obtain the ideal $J$ of isolated roots. 
\par
\vskip 5 pt
\begingroup
\tteight
\baselineskip=10.01pt
\lineskip=0pt
\obeyspaces
i11\ :\ J\ =\ saturate(I,\ ideal(y))\leavevmode\hss\endgraf
\penalty-500\leavevmode\hss\endgraf
\ \ \ \ \ \ \ \ \ \ \ \ \ \ 4\ \ \ \ 3\ \ \ \ \ 2\leavevmode\hss\endgraf
o11\ =\ ideal\ (x\ \ +\ x\ \ +\ 3x\ \ +\ 3x,\ y\ -\ 2x\ -\ 1)\leavevmode\hss\endgraf
\penalty-500\leavevmode\hss\endgraf
o11\ :\ Ideal\ of\ R\leavevmode\hss\endgraf
\endgroup
\penalty-1000
\par
\vskip 1 pt
\noindent
The first polynomial factors completely in ${\mathbb F}_7[x]$,
\par
\vskip 5 pt
\begingroup
\tteight
\baselineskip=10.01pt
\lineskip=0pt
\obeyspaces
i12\ :\ factor(J{\char`\_}0)\leavevmode\hss\endgraf
\penalty-500\leavevmode\hss\endgraf
o12\ =\ (x\ -\ 2)(x\ +\ 1)(x\ +\ 2)(x)(1)\leavevmode\hss\endgraf
\penalty-500\leavevmode\hss\endgraf
o12\ :\ Product\leavevmode\hss\endgraf
\endgroup
\penalty-1000
\par
\vskip 1 pt
\noindent
and so the isolated roots of $I$ are $(0,1),(2,5),(5,4)$, and $(6,6)$. 
\end{example}

Here, the extraneous roots came from a common factor in both
equations.
A less trivial example of this phenomenon will be seen in
Section~\ref{sec:tangent_lines}. 

\subsection{Elimination, rationality, and solving}
Elimination theory\index{elimination theory} can be used to study the
roots of a zero-dimensional ideal $I\subset k[X]$\index{solving polynomial
equations!via elimination}.
A polynomial $h\in k[X]$ defines a map 
$k[y]\rightarrow k[X]$ (by $y\mapsto h$) and a corresponding projection 
$h\colon{\rm Spec}\,k[X]\twoheadrightarrow{\mathbb A}^1$.
The generator $g(y)\in k[y]$ of the 
kernel of the map $k[y]\to k[X]/I$ is called an 
{\it eliminant}\index{eliminant}
and it has the property that ${\mathcal V}(g)=h({\mathcal V}(I))$.
When $h$ is a coordinate function $x_i$, we may consider the eliminant to be
in the polynomial ring $k[x_i]$, and we have 
$\langle g(x_i)\rangle=I\cap k[x_i]$.
The most important result concerning eliminants is the Shape
Lemma\index{Shape Lemma}~\cite{SO:BMMT}. 
\medskip

\noindent{\bf Shape Lemma.} 
{\it
Suppose $h$ is a linear polynomial and $g$ is the corresponding eliminant of
a zero-dimensional ideal $I\subset k[X]$ with $\deg(I)=\deg(g)$. 
Then the roots of $I$ are defined in the splitting
field\index{field!splitting} of $g$ and  
$I$ is radical\index{ideal!radical} if and only if $g$ is square-free.

Suppose further that $h=x_1$ so that $g=g(x_1)$.
Then, in the lexicographic term order 
with $x_1<x_2<\cdots<x_n$, $I$ has a Gr\"obner basis\index{Gr\"obner basis} of
the form: 
\begin{equation}\label{triangular}  
    g(x_1),\ \ x_2-g_2(x_1), \ \ \ldots,\ \ x_n-g_n(x_1)\,,
\end{equation}
where $\deg(g)>\deg(g_i)$ for $i=2,\ldots,n$.
}\medskip

When $k$ is infinite and $I$ is radical, an eliminant $g$ given by a generic
linear polynomial $h$ will satisfy $\deg(g)=\deg(I)$.
Enumerative geometry\index{enumerative geometry} counts solutions
when the fixed figures are generic.
We are similarly concerned with the generic situation of 
$\deg(g)=\deg(I)$.
In this case, eliminants provide a useful computational device to study
further questions about the roots of $I$.
For instance, the Shape Lemma holds for the ideal of Example~\ref{ex:two}.
Its eliminant, which is the polynomial {\tt J{\char`\_}0}, factors completely 
over the ground field ${\mathbb F}_7$, so all four solutions are defined
in ${\mathbb F}_7$.
In Section 4.3, we will use eliminants in another way, to show that
an ideal is radical. 

Given a polynomial $h$ in a zero-dimensional ring $k[X]/I$, the
procedure {\tt eliminant(h, k[y])} finds a linear relation modulo $I$ 
among the powers $1, h, h^2, \ldots, h^d$ of $h$ with $d$ minimal
and returns this as a polynomial in $k[y]$.
This procedure is included in the \Mtwo{}\index{\Mtwo{}}\/ package 
{\tt realroots.m2}.
\par
\vskip 5 pt
\begingroup
\tteight
\baselineskip=10.01pt
\lineskip=0pt
\obeyspaces
i13\ :\ load\ "realroots.m2"\leavevmode\hss\endgraf
\endgroup
\penalty-1000
\par
\vskip 1 pt
\noindent
\par
\vskip 5 pt
\begingroup
\tteight
\baselineskip=10.01pt
\lineskip=0pt
\obeyspaces
i14\ :\ code\ eliminant\leavevmode\hss\endgraf
\penalty-500\leavevmode\hss\endgraf
o14\ =\ --\ code\ for\ eliminant:\leavevmode\hss\endgraf
\ \ \ \ \ \ --\ realroots.m2:65-81\leavevmode\hss\endgraf
\ \ \ \ \ \ eliminant\ =\ (h,\ C)\ ->\ (\leavevmode\hss\endgraf
\ \ \ \ \ \ \ \ \ \ \ Z\ :=\ C{\char`\_}0;\leavevmode\hss\endgraf
\ \ \ \ \ \ \ \ \ \ \ A\ :=\ ring\ h;\leavevmode\hss\endgraf
\ \ \ \ \ \ \ \ \ \ \ assert(\ dim\ A\ ==\ 0\ );\leavevmode\hss\endgraf
\ \ \ \ \ \ \ \ \ \ \ F\ :=\ coefficientRing\ A;\leavevmode\hss\endgraf
\ \ \ \ \ \ \ \ \ \ \ assert(\ isField\ F\ );\leavevmode\hss\endgraf
\ \ \ \ \ \ \ \ \ \ \ assert(\ F\ ==\ coefficientRing\ C\ );\leavevmode\hss\endgraf
\ \ \ \ \ \ \ \ \ \ \ B\ :=\ basis\ A;\leavevmode\hss\endgraf
\ \ \ \ \ \ \ \ \ \ \ d\ :=\ numgens\ source\ B;\leavevmode\hss\endgraf
\ \ \ \ \ \ \ \ \ \ \ M\ :=\ fold((M,\ i)\ ->\ M\ ||\ \leavevmode\hss\endgraf
\ \ \ \ \ \ \ \ \ \ \ \ \ \ \ \ \ \ \ \ \ substitute(contract(B,\ h{\char`\^}(i+1)),\ F),\ \leavevmode\hss\endgraf
\ \ \ \ \ \ \ \ \ \ \ \ \ \ \ \ \ \ \ \ \ substitute(contract(B,\ 1{\char`\_}A),\ F),\ \leavevmode\hss\endgraf
\ \ \ \ \ \ \ \ \ \ \ \ \ \ \ \ \ \ \ \ \ flatten\ subsets(d,\ d));\leavevmode\hss\endgraf
\ \ \ \ \ \ \ \ \ \ \ N\ :=\ ((ker\ transpose\ M)).generators;\leavevmode\hss\endgraf
\ \ \ \ \ \ \ \ \ \ \ P\ :=\ matrix\ {\char`\{}toList\ apply(0..d,\ i\ ->\ Z{\char`\^}i){\char`\}}\ *\ N;\leavevmode\hss\endgraf
\ \ \ \ \ \ \ \ \ \ \ \ \ \ \ \ (flatten\ entries(P)){\char`\_}0\leavevmode\hss\endgraf
\ \ \ \ \ \ \ \ \ \ \ )\leavevmode\hss\endgraf
\penalty-500\leavevmode\hss\endgraf
o14\ :\ Net\leavevmode\hss\endgraf
\endgroup
\penalty-1000
\par
\vskip 1 pt
\noindent
Here, {\tt M} is a matrix whose rows are the normal forms of the
powers  $1$, $h$, $h^2$, $\ldots$, $h^d$ of $h$, for $d$ the degree of the ideal.
The columns of the kernel {\tt N} of {\tt transpose M} are a basis of the
linear relations among these powers. 
The matrix {\tt P} converts these relations into polynomials.
Since {\tt N} is in column echelon form, the initial entry of {\tt P} 
is the relation of minimal degree.
(This method is often faster than na\"\i vely computing the kernel of the
map $k[Z]\to A$ given by $Z\mapsto h$, which is implemented by
{\tt  eliminantNaive(h, Z)}.

Suppose we have an eliminant\index{eliminant} $g(x_1)$ of a zero-dimensional 
ideal  $I\subset k[X]$ with $\deg(g)=\deg(I)$, and we have computed the
lexicographic Gr\"obner basis~(\ref{triangular}). 
Then the roots of $I$ are
\begin{equation}\label{tri_roots}
   \{ (\xi_1,g_2(\xi_1), \ldots, g_n(\xi_1))\mid g(\xi_1)=0\}\,.
\end{equation}

Suppose now that $k={\mathbb Q}$ and we seek floating point approximations
for the (complex) roots of $I$.
Following this method, we first compute floating point solutions to
$g(\xi)=0$, which give all the $x_1$-coordinates of the roots of $I$,  and
then use~(\ref{tri_roots}) to find the other coordinates.
The difficulty here is that enough precision may be lost in evaluating
$g_i(\xi_1)$ so that the result is a poor approximation for the other
components $\xi_i$.

\subsection{Solving with linear algebra}
We describe another method based upon numerical linear algebra.
When $I\subset k[X]$ is zero-dimensional, $A=k[X]/I$ is a finite-dimensional 
$k$-vector space, and {\it any} Gr\"obner basis for $I$ gives an efficient
algorithm to compute ring operations using linear algebra.
In particular, multiplication by $h\in A$ is a linear transformation 
$m_h:A\to A$ and the command {\tt regularRep(h)} from 
{\tt realroots.m2} gives the matrix of $m_h$ in
terms of the standard basis of $A$.
\par
\vskip 5 pt
\begingroup
\tteight
\baselineskip=10.01pt
\lineskip=0pt
\obeyspaces
i15\ :\ code\ regularRep\leavevmode\hss\endgraf
\penalty-500\leavevmode\hss\endgraf
o15\ =\ --\ code\ for\ regularRep:\leavevmode\hss\endgraf
\ \ \ \ \ \ --\ realroots.m2:97-102\leavevmode\hss\endgraf
\ \ \ \ \ \ regularRep\ =\ f\ ->\ (\leavevmode\hss\endgraf
\ \ \ \ \ \ \ \ \ \ \ assert(\ dim\ ring\ f\ ==\ 0\ );\leavevmode\hss\endgraf
\ \ \ \ \ \ \ \ \ \ \ b\ :=\ basis\ ring\ f;\leavevmode\hss\endgraf
\ \ \ \ \ \ \ \ \ \ \ k\ :=\ coefficientRing\ ring\ f;\leavevmode\hss\endgraf
\ \ \ \ \ \ \ \ \ \ \ substitute(contract(transpose\ b,\ f*b),\ k)\leavevmode\hss\endgraf
\ \ \ \ \ \ \ \ \ \ \ )\leavevmode\hss\endgraf
\penalty-500\leavevmode\hss\endgraf
o15\ :\ Net\leavevmode\hss\endgraf
\endgroup
\penalty-1000
\par
\vskip 1 pt
\noindent

Since the action of $A$ on itself is faithful, the minimal polynomial of 
$m_h$ is the eliminant\index{eliminant} corresponding to $h$.
The procedure {\tt charPoly(h, Z)} in {\tt realroots.m2}
computes the characteristic polynomial 
$\det(Z\cdot Id - m_h)$ of $h$.
\par
\vskip 5 pt
\begingroup
\tteight
\baselineskip=10.01pt
\lineskip=0pt
\obeyspaces
i16\ :\ code\ charPoly\leavevmode\hss\endgraf
\penalty-500\leavevmode\hss\endgraf
o16\ =\ --\ code\ for\ charPoly:\leavevmode\hss\endgraf
\ \ \ \ \ \ --\ realroots.m2:108-116\leavevmode\hss\endgraf
\ \ \ \ \ \ charPoly\ =\ (h,\ Z)\ ->\ (\leavevmode\hss\endgraf
\ \ \ \ \ \ \ \ \ \ \ A\ :=\ ring\ h;\leavevmode\hss\endgraf
\ \ \ \ \ \ \ \ \ \ \ F\ :=\ coefficientRing\ A;\leavevmode\hss\endgraf
\ \ \ \ \ \ \ \ \ \ \ S\ :=\ F[Z];\leavevmode\hss\endgraf
\ \ \ \ \ \ \ \ \ \ \ Z\ \ =\ value\ Z;\ \ \ \ \ \leavevmode\hss\endgraf
\ \ \ \ \ \ \ \ \ \ \ mh\ :=\ regularRep(h)\ **\ S;\leavevmode\hss\endgraf
\ \ \ \ \ \ \ \ \ \ \ Idz\ :=\ S{\char`\_}0\ *\ id{\char`\_}(S{\char`\^}(numgens\ source\ mh));\leavevmode\hss\endgraf
\ \ \ \ \ \ \ \ \ \ \ det(Idz\ -\ mh)\leavevmode\hss\endgraf
\ \ \ \ \ \ \ \ \ \ \ )\leavevmode\hss\endgraf
\penalty-500\leavevmode\hss\endgraf
o16\ :\ Net\leavevmode\hss\endgraf
\endgroup
\penalty-1000
\par
\vskip 1 pt
\noindent
When this is the minimal polynomial (the situation of the Shape Lemma),
this procedure often computes the eliminant faster than does 
{\tt eliminant}, and for systems of moderate degree, much faster than
na\"\i vely computing the kernel of the map $k[Z]\to A$ given by $Z\mapsto h$.

The eigenvalues and eigenvectors of $m_h$ give another algorithm for finding
the roots of $I$\index{solving polynomial equations!via eigenvectors}.
The engine for this is the following result\index{Stickelberger's Theorem}.
\medskip

\noindent{\bf Stickelberger's Theorem. }
{\it
Let $h\in A$ and $m_h$ be as above.
Then there is a one-to-one correspondence between eigenvectors
${\bf v}_\xi$ of $m_h$ and roots $\xi$ of $I$, the eigenvalue of $m_h$ on 
${\bf v}_\xi$ is the value $h(\xi)$ of $h$ at $\xi$, and the multiplicity
of this eigenvalue (on the eigenvector ${\bf v}_\xi$) is the
multiplicity of the root $\xi$.
}\medskip

Since the linear transformations $m_h$ for $h\in A$ commute, the
eigenvectors ${\bf v}_\xi$ are common to all $m_h$.
Thus we may compute the roots of a zero-dimensional ideal $I\subset k[X]$
by first computing floating-point approximations to the
eigenvectors ${\bf v}_\xi$ of $m_{x_1}$.
Then the root $\xi\ =\ (\xi_1,\ldots,\xi_n)$ of $I$ corresponding to the
eigenvector ${\bf v}_\xi$ has $i$th coordinate satisfying
\begin{equation}\label{eigenv}
   m_{x_i}\cdot {\bf v}_\xi\ =\ \xi_i \cdot {\bf v}_\xi\,.
\end{equation}
An advantage of this method is that we may use structured numerical linear
algebra after the matrices $m_{x_i}$ are precomputed using exact arithmetic. 
(These matrices are typically sparse and have additional structures which may
be exploited.)
Also, the coordinates $\xi_i$ are {\it linear} functions of the floating
point entries of ${\bf v}_\xi$, which affords greater precision than 
the non-linear evaluations $g_i(\xi_1)$ in the method based upon elimination.
While in principle only one of the $\deg(I)$ components of the vectors
in~(\ref{eigenv}) need be computed, averaging the results from all
components can improve precision.

\subsection{Real Roots}
Determining the real roots of a polynomial system is a challenging problem
with real world applications\index{solving polynomial equations!real solutions}.
When the polynomials come from geometry, this is the main problem of
real enumerative geometry\index{enumerative geometry!real}.
Suppose $k\subset{\mathbb R}$ and $I\subset k[X]$ is zero-dimensional.
If $g$ is an eliminant of $k[X]/I$ 
with $\deg(g)=\deg(I)$, then the real roots of
$g$ are in 1-1 correspondence with the real roots of $I$.
Since there are effective methods for counting the real roots of a univariate
polynomial, eliminants give a na\"\i ve, but useful method for determining the
number of real roots to a polynomial system.
(For some applications of this technique in mathematics,
see~\cite{SO:RS98,SO:So_shap-www,SO:So00b}.) 

The classical symbolic method of Sturm, based upon Sturm sequences,  counts
the number of  real roots of a univariate polynomial in an interval.
When applied to an eliminant satisfying the Shape Lemma, this method counts
the number of real roots of the ideal.
This is implemented in \Mtwo{}\index{\Mtwo{}}\/ via the command
{\tt SturmSequence(f)} of {\tt realroots.m2}
\par
\vskip 5 pt
\begingroup
\tteight
\baselineskip=10.01pt
\lineskip=0pt
\obeyspaces
i17\ :\ code\ SturmSequence\leavevmode\hss\endgraf
\penalty-500\leavevmode\hss\endgraf
o17\ =\ --\ code\ for\ SturmSequence:\leavevmode\hss\endgraf
\ \ \ \ \ \ --\ realroots.m2:120-134\leavevmode\hss\endgraf
\ \ \ \ \ \ SturmSequence\ =\ f\ ->\ (\leavevmode\hss\endgraf
\ \ \ \ \ \ \ \ \ \ \ assert(\ isPolynomialRing\ ring\ f\ );\leavevmode\hss\endgraf
\ \ \ \ \ \ \ \ \ \ \ assert(\ numgens\ ring\ f\ ===\ 1\ );\leavevmode\hss\endgraf
\ \ \ \ \ \ \ \ \ \ \ R\ :=\ ring\ f;\leavevmode\hss\endgraf
\ \ \ \ \ \ \ \ \ \ \ assert(\ char\ R\ ==\ 0\ );\leavevmode\hss\endgraf
\ \ \ \ \ \ \ \ \ \ \ x\ :=\ R{\char`\_}0;\leavevmode\hss\endgraf
\ \ \ \ \ \ \ \ \ \ \ n\ :=\ first\ degree\ f;\leavevmode\hss\endgraf
\ \ \ \ \ \ \ \ \ \ \ c\ :=\ new\ MutableList\ from\ toList\ (0\ ..\ n);\leavevmode\hss\endgraf
\ \ \ \ \ \ \ \ \ \ \ if\ n\ >=\ 0\ then\ (\leavevmode\hss\endgraf
\ \ \ \ \ \ \ \ \ \ \ \ \ \ \ \ c{\char`\#}0\ =\ f;\leavevmode\hss\endgraf
\ \ \ \ \ \ \ \ \ \ \ \ \ \ \ \ if\ n\ >=\ 1\ then\ (\leavevmode\hss\endgraf
\ \ \ \ \ \ \ \ \ \ \ \ \ \ \ \ \ \ \ \ \ c{\char`\#}1\ =\ diff(x,f);\leavevmode\hss\endgraf
\ \ \ \ \ \ \ \ \ \ \ \ \ \ \ \ \ \ \ \ \ scan(2\ ..\ n,\ i\ ->\ c{\char`\#}i\ =\ -\ c{\char`\#}(i-2)\ {\char`\%}\ c{\char`\#}(i-1));\leavevmode\hss\endgraf
\ \ \ \ \ \ \ \ \ \ \ \ \ \ \ \ \ \ \ \ \ ));\leavevmode\hss\endgraf
\ \ \ \ \ \ \ \ \ \ \ toList\ c)\leavevmode\hss\endgraf
\penalty-500\leavevmode\hss\endgraf
o17\ :\ Net\leavevmode\hss\endgraf
\endgroup
\penalty-1000
\par
\vskip 1 pt
\noindent
The last few lines of {\tt SturmSequence} construct the Sturm
sequence\index{Sturm sequence} of the univariate argument $f$:
This is $(f_0, f_1, f_2,\ldots)$ where $f_0=f$, $f_1=f'$, and 
for $i>1$, $f_i$ is the normal form reduction of $-f_{i-2}$ modulo
$f_{i-1}$.
Given any real number $x$, the {\it variation} of $f$ at $x$ is the number of
changes in sign of the sequence $(f_0(x), f_1(x), f_2(x),\ldots)$ obtained by
evaluating the Sturm sequence of $f$ at $x$.
Then the number of real roots of $f$ over an interval $[x,y]$ is the
difference of the variation of $f$ at $x$ and at $y$.

The \Mtwo{}\/ commands  {\tt numRealSturm} and 
{\tt numPosRoots} (and also {\tt numNegRoots}) use this method to respectively
compute the total number of real roots and the number of positive roots of 
a univariate polynomial. 
\par
\vskip 5 pt
\begingroup
\tteight
\baselineskip=10.01pt
\lineskip=0pt
\obeyspaces
i18\ :\ code\ numRealSturm\leavevmode\hss\endgraf
\penalty-500\leavevmode\hss\endgraf
o18\ =\ --\ code\ for\ numRealSturm:\leavevmode\hss\endgraf
\ \ \ \ \ \ --\ realroots.m2:161-165\leavevmode\hss\endgraf
\ \ \ \ \ \ numRealSturm\ =\ f\ ->\ (\leavevmode\hss\endgraf
\ \ \ \ \ \ \ \ \ \ \ c\ :=\ SturmSequence\ f;\leavevmode\hss\endgraf
\ \ \ \ \ \ \ \ \ \ \ variations\ (signAtMinusInfinity\ {\char`\\}\ c)\ \leavevmode\hss\endgraf
\ \ \ \ \ \ \ \ \ \ \ \ \ \ \ -\ variations\ (signAtInfinity\ {\char`\\}\ c)\leavevmode\hss\endgraf
\ \ \ \ \ \ \ \ \ \ \ )\leavevmode\hss\endgraf
\penalty-500\leavevmode\hss\endgraf
o18\ :\ Net\leavevmode\hss\endgraf
\endgroup
\penalty-1000
\par
\vskip 1 pt
\noindent
\par
\vskip 5 pt
\begingroup
\tteight
\baselineskip=10.01pt
\lineskip=0pt
\obeyspaces
i19\ :\ code\ numPosRoots\leavevmode\hss\endgraf
\penalty-500\leavevmode\hss\endgraf
o19\ =\ --\ code\ for\ numPosRoots:\leavevmode\hss\endgraf
\ \ \ \ \ \ --\ realroots.m2:170-174\leavevmode\hss\endgraf
\ \ \ \ \ \ numPosRoots\ =\ f\ ->\ (\ \ \leavevmode\hss\endgraf
\ \ \ \ \ \ \ \ \ \ \ c\ :=\ SturmSequence\ f;\leavevmode\hss\endgraf
\ \ \ \ \ \ \ \ \ \ \ variations\ (signAtZero\ {\char`\\}\ c)\ \leavevmode\hss\endgraf
\ \ \ \ \ \ \ \ \ \ \ \ \ \ \ -\ variations\ (signAtInfinity\ {\char`\\}\ c)\leavevmode\hss\endgraf
\ \ \ \ \ \ \ \ \ \ \ )\leavevmode\hss\endgraf
\penalty-500\leavevmode\hss\endgraf
o19\ :\ Net\leavevmode\hss\endgraf
\endgroup
\penalty-1000
\par
\vskip 1 pt
\noindent
These use the commands {\tt signAt}$*${\tt (f)}, which 
give the sign of ${\tt f}$
at $*$.
(Here, $*$ is one of  {\tt Infinity}, {\tt zero}, or {\tt MinusInfinity}.
Also {\tt variations(c)} computes the 
number of sign changes in the sequence {\tt c}.
\par
\vskip 5 pt
\begingroup
\tteight
\baselineskip=10.01pt
\lineskip=0pt
\obeyspaces
i20\ :\ code\ variations\leavevmode\hss\endgraf
\penalty-500\leavevmode\hss\endgraf
o20\ =\ --\ code\ for\ variations:\leavevmode\hss\endgraf
\ \ \ \ \ \ --\ realroots.m2:187-195\leavevmode\hss\endgraf
\ \ \ \ \ \ variations\ =\ c\ ->\ (\leavevmode\hss\endgraf
\ \ \ \ \ \ \ \ \ \ \ n\ :=\ 0;\leavevmode\hss\endgraf
\ \ \ \ \ \ \ \ \ \ \ last\ :=\ 0;\leavevmode\hss\endgraf
\ \ \ \ \ \ \ \ \ \ \ scan(c,\ x\ ->\ if\ x\ =!=\ 0\ then\ (\leavevmode\hss\endgraf
\ \ \ \ \ \ \ \ \ \ \ \ \ \ \ \ \ \ \ \ \ if\ last\ <\ 0\ and\ x\ >\ 0\ or\ last\ >\ 0\ \leavevmode\hss\endgraf
\ \ \ \ \ \ \ \ \ \ \ \ \ \ \ \ \ \ \ \ \ \ \ \ and\ x\ <\ 0\ then\ n\ =\ n+1;\leavevmode\hss\endgraf
\ \ \ \ \ \ \ \ \ \ \ \ \ \ \ \ \ \ \ \ \ last\ =\ x;\leavevmode\hss\endgraf
\ \ \ \ \ \ \ \ \ \ \ \ \ \ \ \ \ \ \ \ \ ));\leavevmode\hss\endgraf
\ \ \ \ \ \ \ \ \ \ \ n)\leavevmode\hss\endgraf
\penalty-500\leavevmode\hss\endgraf
o20\ :\ Net\leavevmode\hss\endgraf
\endgroup
\penalty-1000
\par
\vskip 1 pt
\noindent

A more sophisticated method to compute the number of real roots which can also
give information about their location uses the rank and
signature\index{bilinear form!signature} of the
symmetric trace form. 
Suppose $I\subset k[X]$ is a zero-dimensional ideal and 
set $A:=k[X]/I$.
For $h\in k[X]$, set $S_h(f,g):={\rm trace}(m_{hfg})$.
It is an easy exercise that $S_h$ is a symmetric bilinear form\index{bilinear
form!symmetric} on $A$.
The procedure {\tt traceForm(h)} in {\tt realroots.m2}
computes this trace form\index{trace form} $S_h$.
\par
\vskip 5 pt
\begingroup
\tteight
\baselineskip=10.01pt
\lineskip=0pt
\obeyspaces
i21\ :\ code\ traceForm\leavevmode\hss\endgraf
\penalty-500\leavevmode\hss\endgraf
o21\ =\ --\ code\ for\ traceForm:\leavevmode\hss\endgraf
\ \ \ \ \ \ --\ realroots.m2:200-208\leavevmode\hss\endgraf
\ \ \ \ \ \ traceForm\ =\ h\ ->\ (\leavevmode\hss\endgraf
\ \ \ \ \ \ \ \ \ \ \ assert(\ dim\ ring\ h\ ==\ 0\ );\leavevmode\hss\endgraf
\ \ \ \ \ \ \ \ \ \ \ b\ \ :=\ basis\ ring\ h;\leavevmode\hss\endgraf
\ \ \ \ \ \ \ \ \ \ \ k\ \ :=\ coefficientRing\ ring\ h;\leavevmode\hss\endgraf
\ \ \ \ \ \ \ \ \ \ \ mm\ :=\ substitute(contract(transpose\ b,\ h\ *\ b\ **\ b),\ k);\leavevmode\hss\endgraf
\ \ \ \ \ \ \ \ \ \ \ tr\ :=\ matrix\ {\char`\{}apply(first\ entries\ b,\ x\ ->\leavevmode\hss\endgraf
\ \ \ \ \ \ \ \ \ \ \ \ \ \ \ \ \ \ \ \ \ trace\ regularRep\ x){\char`\}};\leavevmode\hss\endgraf
\ \ \ \ \ \ \ \ \ \ \ adjoint(tr\ *\ mm,\ source\ tr,\ source\ tr)\leavevmode\hss\endgraf
\ \ \ \ \ \ \ \ \ \ \ )\leavevmode\hss\endgraf
\penalty-500\leavevmode\hss\endgraf
o21\ :\ Net\leavevmode\hss\endgraf
\endgroup
\penalty-1000
\par
\vskip 1 pt
\noindent
The value of this construction is the following theorem.

\begin{theorem}[\cite{SO:BW,SO:PRS}]\label{t:PRS}
Suppose $k\subset{\mathbb R}$ and $I$ is a zero-dimensional ideal in
$k[x_1,\ldots,x_n]$ and consider ${\mathcal V}(I)\subset {\mathbb C}^n$. 
Then, for $h\in k[x_1,\ldots,x_n]$, the signature $\sigma(S_h)$ and rank 
$\rho(S_h)$ of the bilinear form $S_h$ satisfy
\begin{eqnarray*}
\sigma(S_h)&=&\#\{a\in{\mathcal V}(I)\cap{\mathbb R}^n:h(a)>0\}
            - \#\{a\in{\mathcal V}(I)\cap{\mathbb R}^n:h(a)<0\}\,\\
\rho(S_h)&=&\#\{a\in{\mathcal V}(I):h(a)\neq0\}\,.
\end{eqnarray*}
\end{theorem}

That is, the rank of $S_h$ counts roots in 
${\mathbb C}^n-{\mathcal V}(h)$, and its signature counts the real roots
weighted by the sign of $h$ (which is $-1$, $0$, or $1$) at each root.
The command {\tt traceFormSignature(h)} in {\tt realroots.m2} returns the
rank and  signature of the trace form $S_h$.
\par
\vskip 5 pt
\begingroup
\tteight
\baselineskip=10.01pt
\lineskip=0pt
\obeyspaces
i22\ :\ code\ traceFormSignature\leavevmode\hss\endgraf
\penalty-500\leavevmode\hss\endgraf
o22\ =\ --\ code\ for\ traceFormSignature:\leavevmode\hss\endgraf
\ \ \ \ \ \ --\ realroots.m2:213-224\leavevmode\hss\endgraf
\ \ \ \ \ \ traceFormSignature\ =\ h\ ->\ (\leavevmode\hss\endgraf
\ \ \ \ \ \ \ \ \ \ \ A\ :=\ ring\ h;\leavevmode\hss\endgraf
\ \ \ \ \ \ \ \ \ \ \ assert(\ dim\ A\ ==\ 0\ );\leavevmode\hss\endgraf
\ \ \ \ \ \ \ \ \ \ \ assert(\ char\ A\ ==\ 0\ );\leavevmode\hss\endgraf
\ \ \ \ \ \ \ \ \ \ \ S\ :=\ QQ[Z];\leavevmode\hss\endgraf
\ \ \ \ \ \ \ \ \ \ \ TrF\ :=\ traceForm(h)\ **\ S;\leavevmode\hss\endgraf
\ \ \ \ \ \ \ \ \ \ \ IdZ\ :=\ Z\ *\ id{\char`\_}(S{\char`\^}(numgens\ source\ TrF));\leavevmode\hss\endgraf
\ \ \ \ \ \ \ \ \ \ \ f\ :=\ det(TrF\ -\ IdZ);\leavevmode\hss\endgraf
\ \ \ \ \ \ \ \ \ \ \ <<\ "The\ trace\ form\ S{\char`\_}h\ with\ h\ =\ "\ <<\ h\ <<\ \leavevmode\hss\endgraf
\ \ \ \ \ \ \ \ \ \ \ \ \ "\ has\ rank\ "\ <<\ rank(TrF)\ <<\ "\ and\ signature\ "\ <<\ \leavevmode\hss\endgraf
\ \ \ \ \ \ \ \ \ \ \ \ \ numPosRoots(f)\ -\ numNegRoots(f)\ <<\ endl\leavevmode\hss\endgraf
\ \ \ \ \ \ \ \ \ \ \ )\leavevmode\hss\endgraf
\penalty-500\leavevmode\hss\endgraf
o22\ :\ Net\leavevmode\hss\endgraf
\endgroup
\penalty-1000
\par
\vskip 1 pt
\noindent
The \Mtwo{}\/ command {\tt numRealTrace(A)} simply returns the number of
real roots of $I$, given ${\tt A}=k[X]/I$.  
\par
\vskip 5 pt
\begingroup
\tteight
\baselineskip=10.01pt
\lineskip=0pt
\obeyspaces
i23\ :\ code\ numRealTrace\leavevmode\hss\endgraf
\penalty-500\leavevmode\hss\endgraf
o23\ =\ --\ code\ for\ numRealTrace:\leavevmode\hss\endgraf
\ \ \ \ \ \ --\ realroots.m2:229-237\leavevmode\hss\endgraf
\ \ \ \ \ \ numRealTrace\ =\ A\ ->\ (\leavevmode\hss\endgraf
\ \ \ \ \ \ \ \ \ \ \ assert(\ dim\ A\ ==\ 0\ );\leavevmode\hss\endgraf
\ \ \ \ \ \ \ \ \ \ \ assert(\ char\ A\ ==\ 0\ );\leavevmode\hss\endgraf
\ \ \ \ \ \ \ \ \ \ \ S\ :=\ QQ[Z];\leavevmode\hss\endgraf
\ \ \ \ \ \ \ \ \ \ \ TrF\ :=\ traceForm(1{\char`\_}A)\ **\ S;\leavevmode\hss\endgraf
\ \ \ \ \ \ \ \ \ \ \ IdZ\ :=\ Z\ *\ id{\char`\_}(S{\char`\^}(numgens\ source\ TrF));\leavevmode\hss\endgraf
\ \ \ \ \ \ \ \ \ \ \ f\ :=\ det(TrF\ -\ IdZ);\leavevmode\hss\endgraf
\ \ \ \ \ \ \ \ \ \ \ numPosRoots(f)-numNegRoots(f)\leavevmode\hss\endgraf
\ \ \ \ \ \ \ \ \ \ \ )\leavevmode\hss\endgraf
\penalty-500\leavevmode\hss\endgraf
o23\ :\ Net\leavevmode\hss\endgraf
\endgroup
\penalty-1000
\par
\vskip 1 pt
\noindent

\begin{example}
We illustrate these methods on the following polynomial system.
\par
\vskip 5 pt
\begingroup
\tteight
\baselineskip=10.01pt
\lineskip=0pt
\obeyspaces
i24\ :\ R\ =\ QQ[x,\ y];\leavevmode\hss\endgraf
\endgroup
\penalty-1000
\par
\vskip 1 pt
\noindent
\par
\vskip 5 pt
\begingroup
\tteight
\baselineskip=10.01pt
\lineskip=0pt
\obeyspaces
i25\ :\ I\ =\ ideal\ (1\ -\ x{\char`\^}2*y\ +\ 2*x*y{\char`\^}2,\ \ y\ -\ 2*x\ -\ x*y\ +\ x{\char`\^}2);\leavevmode\hss\endgraf
\penalty-500\leavevmode\hss\endgraf
o25\ :\ Ideal\ of\ R\leavevmode\hss\endgraf
\endgroup
\penalty-1000
\par
\vskip 1 pt
\noindent
The ideal $I$ has dimension zero and degree 5.
\par
\vskip 5 pt
\begingroup
\tteight
\baselineskip=10.01pt
\lineskip=0pt
\obeyspaces
i26\ :\ dim\ I,\ degree\ I\leavevmode\hss\endgraf
\penalty-500\leavevmode\hss\endgraf
o26\ =\ (0,\ 5)\leavevmode\hss\endgraf
\penalty-500\leavevmode\hss\endgraf
o26\ :\ Sequence\leavevmode\hss\endgraf
\endgroup
\penalty-1000
\par
\vskip 1 pt
\noindent
We compare the two methods to compute the eliminant of $x$ in 
the ring $R/I$.
\par
\vskip 5 pt
\begingroup
\tteight
\baselineskip=10.01pt
\lineskip=0pt
\obeyspaces
i27\ :\ A\ =\ R/I;\leavevmode\hss\endgraf
\endgroup
\penalty-1000
\par
\vskip 1 pt
\noindent
\par
\vskip 5 pt
\begingroup
\tteight
\baselineskip=10.01pt
\lineskip=0pt
\obeyspaces
i28\ :\ time\ g\ =\ eliminant(x,\ QQ[Z])\leavevmode\hss\endgraf
\ \ \ \ \ --\ used\ 0.03\ seconds\leavevmode\hss\endgraf
\penalty-500\leavevmode\hss\endgraf
\ \ \ \ \ \ \ 5\ \ \ \ \ 4\ \ \ \ \ 3\ \ \ \ 2\leavevmode\hss\endgraf
o28\ =\ Z\ \ -\ 5Z\ \ +\ 6Z\ \ +\ Z\ \ -\ 2Z\ +\ 1\leavevmode\hss\endgraf
\penalty-500\leavevmode\hss\endgraf
o28\ :\ QQ\ [Z]\leavevmode\hss\endgraf
\endgroup
\penalty-1000
\par
\vskip 1 pt
\noindent
\par
\vskip 5 pt
\begingroup
\tteight
\baselineskip=10.01pt
\lineskip=0pt
\obeyspaces
i29\ :\ time\ g\ =\ charPoly(x,\ Z)\leavevmode\hss\endgraf
\ \ \ \ \ --\ used\ 0.01\ seconds\leavevmode\hss\endgraf
\penalty-500\leavevmode\hss\endgraf
\ \ \ \ \ \ \ 5\ \ \ \ \ 4\ \ \ \ \ 3\ \ \ \ 2\leavevmode\hss\endgraf
o29\ =\ Z\ \ -\ 5Z\ \ +\ 6Z\ \ +\ Z\ \ -\ 2Z\ +\ 1\leavevmode\hss\endgraf
\penalty-500\leavevmode\hss\endgraf
o29\ :\ QQ\ [Z]\leavevmode\hss\endgraf
\endgroup
\penalty-1000
\par
\vskip 1 pt
\noindent
The eliminant has 3 real roots, which we test in two different ways.
\par
\vskip 5 pt
\begingroup
\tteight
\baselineskip=10.01pt
\lineskip=0pt
\obeyspaces
i30\ :\ numRealSturm(g),\ numRealTrace(A)\leavevmode\hss\endgraf
\penalty-500\leavevmode\hss\endgraf
o30\ =\ (3,\ 3)\leavevmode\hss\endgraf
\penalty-500\leavevmode\hss\endgraf
o30\ :\ Sequence\leavevmode\hss\endgraf
\endgroup
\penalty-1000
\par
\vskip 1 pt
\noindent
We use Theorem~\ref{t:PRS} to isolate these roots in the $x,y$-plane.
\par
\vskip 5 pt
\begingroup
\tteight
\baselineskip=10.01pt
\lineskip=0pt
\obeyspaces
i31\ :\ traceFormSignature(x*y);\leavevmode\hss\endgraf
The\ trace\ form\ S{\char`\_}h\ with\ h\ =\ x*y\ has\ rank\ 5\ and\ signature\ 3\leavevmode\hss\endgraf
\endgroup
\penalty-1000
\par
\vskip 1 pt
\noindent
Thus all 3 real roots lie in the first and third
quadrants (where $xy>0$).
We isolate these further.
\par
\vskip 5 pt
\begingroup
\tteight
\baselineskip=10.01pt
\lineskip=0pt
\obeyspaces
i32\ :\ traceFormSignature(x\ -\ 2);\leavevmode\hss\endgraf
The\ trace\ form\ S{\char`\_}h\ with\ h\ =\ x\ -\ 2\ has\ rank\ 5\ and\ signature\ 1\leavevmode\hss\endgraf
\endgroup
\penalty-1000
\par
\vskip 1 pt
\noindent
This shows that two roots lie in the first quadrant with $x>2$ and one lies
in the third.
Finally, one of the roots lies in the triangle $y>0$, $x>2$, and $x+y<3$.
\par
\vskip 5 pt
\begingroup
\tteight
\baselineskip=10.01pt
\lineskip=0pt
\obeyspaces
i33\ :\ traceFormSignature(x\ +\ y\ -\ 3);\leavevmode\hss\endgraf
The\ trace\ form\ S{\char`\_}h\ with\ h\ =\ x\ +\ y\ -\ 3\ has\ rank\ 5\ and\ signature\ -1\leavevmode\hss\endgraf
\endgroup
\penalty-1000
\par
\vskip 1 pt
\noindent

Figure~\ref{fig:roots} shows these three roots (dots), as well as the
lines $x+y=3$ and $x=2$.
\begin{figure}[htb]
$$
  %
  %
  \setlength{\unitlength}{0.96pt}
  \begin{picture}(220,110)(-60,-50)
   \put(0,-50){\vector(0,1){110}} \put(0,0){\vector(0,-1){50}}
   \put(-11,47){$y$}
   \put(-5,-40){\line(1,0){10}} \put(6,-42.8){$-1$}
   \put(-5,40){\line(1,0){10}}  \put(6, 36.5){$1$}
   \put(-60,0){\vector(1,0){220}} \put(0,0){\vector(-1,0){60}}
   \put(-40,-5){\line(0,1){10}} \put(-50,-15){$-1$}
   \put(40,-5){\line(0,1){10}}  \put(37.5,-15){$1$}
   \put(120,-5){\line(0,1){10}} \put(117,-15){$3$}
   \put(-51,5){$x$}
  
  \thicklines
  \put(80,-50){\line(0,1){110}}  \put(85,-35){$x=2$}
  \put(160,-40){\line(-1,1){100}}
  \put(100,25){$x+y=3$}
  
  \put(-26.2, -42.1){\circle*{2}}
  \put(86.3, 11.8){\circle*{2}}
  \put(112.8, 50.8){\circle*{2}}
  \end{picture}
$$
\caption{Location of roots\label{fig:roots}}
\end{figure}
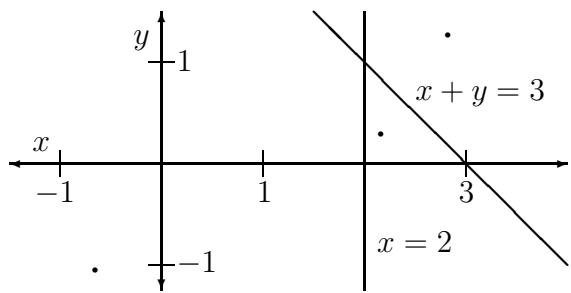
\end{example}

\subsection{Homotopy methods}
We describe symbolic-numeric 
{\it homotopy continuation methods}\index{homotopy continuation} 
for finding approximate complex solutions to a system of
equations\index{solving polynomial equations!via numerical homotopy}.
These exploit the traditional principles of conservation of number and 
specialization from enumerative geometry\index{enumerative geometry}.

Suppose we seek the isolated solutions of a system $F(X)=0$
where $F=(f_1,\ldots,f_n)$ are polynomials in the variables
$X=(x_1,\ldots,x_N)$.
First, a {\em homotopy} $H(X,t)$ is found with the following properties:
\begin{enumerate}
  \item $H(X,1)= F(X)$. 
  \item The isolated solutions of the {\it start system} $H(X,0)=0$ are known.
  \item The system $H(X,t)=0$ defines finitely many (complex) curves, 
        and each isolated solution of the original system $F(X)=0$ is
        connected to an isolated solution $\sigma_i(0)$ of $H(X,0)=0$ along
        one of these curves. 
\end{enumerate}
Next, choose a generic smooth path $\gamma(t)$ from 0 to 1 in the complex
plane.
Lifting $\gamma$ to the curves $H(X,t)=0$ gives 
smooth paths $\sigma_i(t)$ connecting each solution
$\sigma_i(0)$ of the start system to a solution of the original system.
The path $\gamma$ must avoid the finitely many points in ${\mathbb C}$ over
which the curves are singular or meet other components of the solution set
$H(X,t)=0$.

Numerical path continuation is used to trace each path
$\sigma_i(t)$ from $t=0$ to $t=1$.
When there are fewer solutions to $F(X)=0$ than to 
$H(X,0)=0$, some paths will diverge or become singular as
$t\rightarrow 1$, and it is expensive to trace such a path.
The homotopy is {\it optimal}\index{homotopy!optimal} when this does not
occur. 

When $N=n$ and the $f_i$ are generic, set
$G(X):=(g_1,\ldots,g_n)$ with $g_i=(x_i-1)(x_i-2)\cdots(x_i-d_i)$
where $d_i:=\deg(f_i)$.
Then the {\it B\'ezout homotopy}\index{homotopy!B\'ezout} 
$$
  H(X,t)\quad :=\quad tF(X)\ +\ 
  (1-t)G(X)
$$
is optimal.
This homotopy furnishes an effective demonstration of
the bound in B\'ezout's Theorem\index{B\'ezout Theorem} for the number of
solutions to $F(X)=0$.

When the polynomial system is deficient, the B\'ezout homotopy is not optimal.
When $n>N$ (often the case in geometric examples), 
the B\'ezout homotopy  does not apply.
In either case, a different strategy is needed.
Present optimal homotopies for such systems all exploit some structure of
the systems they are designed to solve.
The current state-of-the-art is described in~\cite{SO:Ver99}.

\begin{example}\label{example:Groebner}
The Gr\"obner homotopy\index{homotopy!Gr\"obner}~\cite{SO:HSS} is an optimal 
homotopy\index{homotopy!optimal} that exploits a square-free initial
ideal\index{initial ideal!square-free}.
Suppose our system has the form
$$
  F\ :=\ g_1(X),\ldots,g_m(X),\ \Lambda_1(X),\ldots,\Lambda_d(X)
$$
where $g_1(X),\ldots,g_m(X)$ form a Gr\"obner basis for an ideal $I$ 
with respect to a given term order $\prec$, $\Lambda_1,\ldots,\Lambda_d$ are
linear forms with $d=\dim({\mathcal V}(I))$, {\it and}\/ we assume that 
the initial ideal ${\rm in}_\prec I$ is square-free.
This last, restrictive, hypothesis occurs for
certain determinantal varieties.

As in~\cite[Chapter 15]{SO:MR97a:13001}, there exist polynomials
$g_i(X,t)$ interpolating between $g_i(X)$ and their initial terms
${\rm in}_\prec g_i(X)$
$$
  g_i(X;1)\ =\ g_i(X) \qquad\mbox{and}\qquad
  g_i(X;0) \ =\ {\rm in}_\prec g_i(X)
$$
so that $\langle g_1(X,t),\ldots,g_m(X,t)\rangle$ is a flat family
with generic fibre isomorphic to $I$ and special fibre 
${\rm in}_\prec I$.
The {\it Gr\"obner homotopy} is
$$
  H(X,t)\ :=\ 
  g_1(X,t),\ldots,g_m(X,t),\ \Lambda_1(X),\ldots,\Lambda_d(X).
$$
Since ${\rm in}_\prec I$ is square-free, 
${\mathcal V}({\rm in}_\prec I)$ is a union of 
$\deg(I)$-many coordinate $d$-planes.
We solve the start system by linear algebra.
This conceptually simple homotopy is in general not
efficient as it is typically overdetermined.
\end{example}

\section{Some enumerative geometry}\label{sec:enumerative}

We use the tools we have developed to explore the enumerative geometric
problems of cylinders meeting 5 general points and lines tangent to
4 spheres\index{enumerative geometry}\index{enumerative problem}. 

\subsection{Cylinders meeting 5 points}\label{sec:cylinder}
A {\it cylinder}\index{cylinder} is the locus of points equidistant from a
fixed line in ${\mathbb R}^3$. 
The Grassmannian\index{Grassmannian} of lines in 3-space is 4-dimensional,
which implies that 
the space of cylinders is 5-dimensional, and so we expect that 5 points in  
${\mathbb R}^3$ will determine finitely many cylinders.
That is, there should be finitely many lines equidistant from 5 general points.
The question is: How many cylinders/lines, and how many of them can be real?

Bottema\index{Bottema, O.} and Veldkamp\index{Veldkamp, G.}~\cite{SO:BV77}
show there are 6 {\it complex} cylinders 
and Lichtblau\index{Lichtblau, D.}~\cite{SO:Li00} observes that if the 5
points are the vertices of 
a bipyramid consisting of 2 regular tetrahedra sharing a common face, then
all 6 will be real.
We check this reality on a configuration with less symmetry (so the Shape
Lemma holds).

If the axial line has direction ${\bf V}$ and contains the point ${\bf P}$
(and hence has parameterization ${\bf P}+t{\bf V}$), and if $r$ is the squared
radius, then the cylinder\index{cylinder} is the set of points ${\bf X}$
satisfying 
$$
  0 \ =\ r - 
  \left\| {\bf X} - {\bf P} - \frac{{\bf V}\cdot({\bf X} - {\bf P})}%
     {\|{\bf V}\|^2}\,{\bf V} \right\|^2\ .
$$
Expanding and clearing the denominator of $\|{\bf V}\|^2$ yields
\begin{equation}\label{eq:cylinder}
  0 \ =\ r \|{\bf V}\|^2 + 
         [{\bf V}\cdot({\bf X} - {\bf P})]^2 - 
        \|{\bf X} - {\bf P}\|^2\, \|{\bf V}\|^2\,.
\end{equation}
We consider cylinders containing the following 5 points, which form an
asymmetric bipyramid.
\par
\vskip 5 pt
\begingroup
\tteight
\baselineskip=10.01pt
\lineskip=0pt
\obeyspaces
i34\ :\ Points\ =\ {\char`\{}{\char`\{}2,\ 2,\ \ 0\ {\char`\}},\ {\char`\{}1,\ -2,\ \ 0{\char`\}},\ {\char`\{}-3,\ 0,\ 0{\char`\}},\ \leavevmode\hss\endgraf
\ \ \ \ \ \ \ \ \ \ \ \ \ \ \ \ {\char`\{}0,\ 0,\ 5/2{\char`\}},\ {\char`\{}0,\ \ 0,\ -3{\char`\}}{\char`\}};\leavevmode\hss\endgraf
\endgroup
\penalty-1000
\par
\vskip 1 pt
\noindent
Suppose that ${\bf P}=(0,y_{11},y_{12})$ and ${\bf V}=(1,y_{21},y_{22})$.
\par
\vskip 5 pt
\begingroup
\tteight
\baselineskip=10.01pt
\lineskip=0pt
\obeyspaces
i35\ :\ R\ =\ QQ[r,\ y11,\ y12,\ y21,\ y22];\leavevmode\hss\endgraf
\endgroup
\penalty-1000
\par
\vskip 1 pt
\noindent
\par
\vskip 5 pt
\begingroup
\tteight
\baselineskip=10.01pt
\lineskip=0pt
\obeyspaces
i36\ :\ P\ =\ matrix{\char`\{}{\char`\{}0,\ y11,\ y12{\char`\}}{\char`\}};\leavevmode\hss\endgraf
\penalty-500\leavevmode\hss\endgraf
\ \ \ \ \ \ \ \ \ \ \ \ \ \ 1\ \ \ \ \ \ \ 3\leavevmode\hss\endgraf
o36\ :\ Matrix\ R\ \ <---\ R\leavevmode\hss\endgraf
\endgroup
\penalty-1000
\par
\vskip 1 pt
\noindent
\par
\vskip 5 pt
\begingroup
\tteight
\baselineskip=10.01pt
\lineskip=0pt
\obeyspaces
i37\ :\ V\ =\ matrix{\char`\{}{\char`\{}1,\ y21,\ y22{\char`\}}{\char`\}};\leavevmode\hss\endgraf
\penalty-500\leavevmode\hss\endgraf
\ \ \ \ \ \ \ \ \ \ \ \ \ \ 1\ \ \ \ \ \ \ 3\leavevmode\hss\endgraf
o37\ :\ Matrix\ R\ \ <---\ R\leavevmode\hss\endgraf
\endgroup
\penalty-1000
\par
\vskip 1 pt
\noindent
We construct the ideal given by evaluating the
polynomial~(\ref{eq:cylinder}) at each of the five points.
\par
\vskip 5 pt
\begingroup
\tteight
\baselineskip=10.01pt
\lineskip=0pt
\obeyspaces
i38\ :\ Points\ =\ matrix\ Points\ **\ R;\leavevmode\hss\endgraf
\penalty-500\leavevmode\hss\endgraf
\ \ \ \ \ \ \ \ \ \ \ \ \ \ 5\ \ \ \ \ \ \ 3\leavevmode\hss\endgraf
o38\ :\ Matrix\ R\ \ <---\ R\leavevmode\hss\endgraf
\endgroup
\penalty-1000
\par
\vskip 1 pt
\noindent
\par
\vskip 5 pt
\begingroup
\tteight
\baselineskip=10.01pt
\lineskip=0pt
\obeyspaces
i39\ :\ I\ =\ ideal\ apply(0..4,\ i\ ->\ (\leavevmode\hss\endgraf
\ \ \ \ \ \ \ \ \ \ \ \ \ \ \ \ X\ :=\ Points{\char`\^}{\char`\{}i{\char`\}};\leavevmode\hss\endgraf
\ \ \ \ \ \ \ \ \ \ \ \ \ \ \ \ r\ *\ (V\ *\ transpose\ V)\ \ +\leavevmode\hss\endgraf
\ \ \ \ \ \ \ \ \ \ \ \ \ \ \ \ \ ((X\ -\ P)\ *\ transpose\ V){\char`\^}2)\ -\leavevmode\hss\endgraf
\ \ \ \ \ \ \ \ \ \ \ \ \ \ \ \ \ ((X\ -\ P)\ *\ transpose(X\ -\ P))\ *\ (V\ *\ transpose\ V)\leavevmode\hss\endgraf
\ \ \ \ \ \ \ \ \ \ \ \ \ \ \ \ );\leavevmode\hss\endgraf
\penalty-500\leavevmode\hss\endgraf
o39\ :\ Ideal\ of\ R\leavevmode\hss\endgraf
\endgroup
\penalty-1000
\par
\vskip 1 pt
\noindent
This ideal has dimension 0 and degree 6.
\par
\vskip 5 pt
\begingroup
\tteight
\baselineskip=10.01pt
\lineskip=0pt
\obeyspaces
i40\ :\ dim\ I,\ degree\ I\leavevmode\hss\endgraf
\penalty-500\leavevmode\hss\endgraf
o40\ =\ (0,\ 6)\leavevmode\hss\endgraf
\penalty-500\leavevmode\hss\endgraf
o40\ :\ Sequence\leavevmode\hss\endgraf
\endgroup
\penalty-1000
\par
\vskip 1 pt
\noindent
There are 6 real roots, and they correspond to real cylinders (with $r>0$).
\par
\vskip 5 pt
\begingroup
\tteight
\baselineskip=10.01pt
\lineskip=0pt
\obeyspaces
i41\ :\ A\ =\ R/I;\ numPosRoots(charPoly(r,\ Z))\leavevmode\hss\endgraf
\penalty-500\leavevmode\hss\endgraf
o42\ =\ 3\leavevmode\hss\endgraf
\endgroup
\penalty-1000
\par
\vskip 1 pt
\noindent

\subsection{Lines tangent to 4 spheres}\label{sec:12lines}
We now ask for the lines having a fixed distance from 4 general points.
Equivalently, these are the lines mutually tangent to 4 spheres\index{sphere}.
Since the Grassmannian\index{Grassmannian} of lines is four-dimensional, we
expect there to be only finitely many such lines.
Macdonald\index{Macdonald, I.}, Pach\index{Pach, J.}, and
Theobald\index{Theobald, Th.}~\cite{SO:MPT00} show that there
are indeed 12 lines, and that all 12 may be real.
This problem makes geometric sense over any field $k$ not of characteristic
2, and the derivation of the number 12 is also valid for algebraically
closed\index{field!algebraically closed} 
fields not of characteristic 2.

A sphere in $k^3$ is ${\mathcal V}(q(1,{\bf x}))$, where $q$ is a
quadratic form\index{quadratic form} on $k^4$ and ${\bf x}\in k^3$.
If our field does not have characteristic 2, then there 
is a symmetric $4\times 4$ matrix $M$ such that 
$q({\bf u})={\bf u}M{\bf u}^t$.

A line $\ell$ having direction ${\bf V}$ and containing the point ${\bf P}$ 
is tangent to the sphere defined by $q$ when the univariate polynomial in $s$ 
$$
  q( (1,{\bf P})+s(0,{\bf V}) )\ =\ 
  q(1,{\bf P}) + 2s (1,{\bf P})M (0,{\bf V})^t + s^2q(0,{\bf V})\,,
$$
has a double root.
Thus its discriminant\index{discriminant} vanishes, giving the equation
\begin{equation}\label{eq:sphere}
  \left( (1,{\bf P})M(0,{\bf V})^t\right)^2 \ -\ 
   (1,{\bf P})M (1,{\bf P})^t\cdot(0,{\bf V})M (0,{\bf V})^t 
     \ =\ 0\,.
\end{equation}

The matrix $M$ of the quadratic form $q$ of the sphere with
center $(a,b,c)$ and squared radius $r$ is constructed by 
{\tt Sphere(a,b,c,r)}.
\par
\vskip 5 pt
\begingroup
\tteight
\baselineskip=10.01pt
\lineskip=0pt
\obeyspaces
i43\ :\ Sphere\ =\ (a,\ b,\ c,\ r)\ ->\ (\leavevmode\hss\endgraf
\ \ \ \ \ \ \ \ \ \ \ \ \ \ matrix{\char`\{}{\char`\{}a{\char`\^}2\ +\ b{\char`\^}2\ +\ c{\char`\^}2\ -\ r\ ,-a\ ,-b\ ,-c\ {\char`\}},\leavevmode\hss\endgraf
\ \ \ \ \ \ \ \ \ \ \ \ \ \ \ \ \ \ \ \ \ {\char`\{}\ \ \ \ \ \ \ \ \ -a\ \ \ \ \ \ \ \ \ ,\ 1\ ,\ 0\ ,\ 0\ {\char`\}},\leavevmode\hss\endgraf
\ \ \ \ \ \ \ \ \ \ \ \ \ \ \ \ \ \ \ \ \ {\char`\{}\ \ \ \ \ \ \ \ \ -b\ \ \ \ \ \ \ \ \ ,\ 0\ ,\ 1\ ,\ 0\ {\char`\}},\leavevmode\hss\endgraf
\ \ \ \ \ \ \ \ \ \ \ \ \ \ \ \ \ \ \ \ \ {\char`\{}\ \ \ \ \ \ \ \ \ -c\ \ \ \ \ \ \ \ \ ,\ 0\ ,\ 0\ ,\ 1\ {\char`\}}{\char`\}}\leavevmode\hss\endgraf
\ \ \ \ \ \ \ \ \ \ \ \ \ \ );\leavevmode\hss\endgraf
\endgroup
\penalty-1000
\par
\vskip 1 pt
\noindent
If a line $\ell$ contains the point ${\bf P}=(0,y_{11},y_{12})$ 
and $\ell$ has direction ${\bf V} = (1,y_{21},y_{22})$, then 
{\tt tangentTo(M)} is the equation for $\ell$ to be tangent to the 
quadric $uMu^T=0$ determined by the matrix $M$.
\par
\vskip 5 pt
\begingroup
\tteight
\baselineskip=10.01pt
\lineskip=0pt
\obeyspaces
i44\ :\ R\ =\ QQ[y11,\ y12,\ y21,\ y22];\leavevmode\hss\endgraf
\endgroup
\penalty-1000
\par
\vskip 1 pt
\noindent
\par
\vskip 5 pt
\begingroup
\tteight
\baselineskip=10.01pt
\lineskip=0pt
\obeyspaces
i45\ :\ tangentTo\ =\ (M)\ ->\ (\leavevmode\hss\endgraf
\ \ \ \ \ \ \ \ \ \ \ P\ :=\ matrix{\char`\{}{\char`\{}1,\ 0,\ y11,\ y12{\char`\}}{\char`\}};\leavevmode\hss\endgraf
\ \ \ \ \ \ \ \ \ \ \ V\ :=\ matrix{\char`\{}{\char`\{}0,\ 1,\ y21,\ y22{\char`\}}{\char`\}};\leavevmode\hss\endgraf
\ \ \ \ \ \ \ \ \ \ \ (P\ *\ M\ *\ transpose\ V){\char`\^}2\ -\ \leavevmode\hss\endgraf
\ \ \ \ \ \ \ \ \ \ \ \ \ (P\ *\ M\ *\ transpose\ P)\ *\ (V\ *\ M\ *\ transpose\ V)\leavevmode\hss\endgraf
\ \ \ \ \ \ \ \ \ \ \ );\leavevmode\hss\endgraf
\endgroup
\penalty-1000
\par
\vskip 1 pt
\noindent
The ideal of lines having distance $\sqrt{5}$ from the four points
$(0,0,0)$, $(4,1,1)$, $(1,4,1)$, and $(1,1,4)$ has dimension zero and degree 12.
\par
\vskip 5 pt
\begingroup
\tteight
\baselineskip=10.01pt
\lineskip=0pt
\obeyspaces
i46\ :\ I\ =\ ideal\ (tangentTo(Sphere(0,0,0,5)),\ \leavevmode\hss\endgraf
\ \ \ \ \ \ \ \ \ \ \ \ \ \ \ \ \ tangentTo(Sphere(4,1,1,5)),\ \leavevmode\hss\endgraf
\ \ \ \ \ \ \ \ \ \ \ \ \ \ \ \ \ tangentTo(Sphere(1,4,1,5)),\ \leavevmode\hss\endgraf
\ \ \ \ \ \ \ \ \ \ \ \ \ \ \ \ \ tangentTo(Sphere(1,1,4,5)));\leavevmode\hss\endgraf
\penalty-500\leavevmode\hss\endgraf
o46\ :\ Ideal\ of\ R\leavevmode\hss\endgraf
\endgroup
\penalty-1000
\par
\vskip 1 pt
\noindent
\par
\vskip 5 pt
\begingroup
\tteight
\baselineskip=10.01pt
\lineskip=0pt
\obeyspaces
i47\ :\ dim\ I,\ degree\ I\leavevmode\hss\endgraf
\penalty-500\leavevmode\hss\endgraf
o47\ =\ (0,\ 12)\leavevmode\hss\endgraf
\penalty-500\leavevmode\hss\endgraf
o47\ :\ Sequence\leavevmode\hss\endgraf
\endgroup
\penalty-1000
\par
\vskip 1 pt
\noindent
Thus there are 12 lines whose distance from those 4 points is $\sqrt{5}$.
We check that all 12 are real.
\par
\vskip 5 pt
\begingroup
\tteight
\baselineskip=10.01pt
\lineskip=0pt
\obeyspaces
i48\ :\ A\ =\ R/I;\leavevmode\hss\endgraf
\endgroup
\penalty-1000
\par
\vskip 1 pt
\noindent
\par
\vskip 5 pt
\begingroup
\tteight
\baselineskip=10.01pt
\lineskip=0pt
\obeyspaces
i49\ :\ numRealSturm(eliminant(y11\ -\ y12\ +\ y21\ +\ y22,\ QQ[Z]))\leavevmode\hss\endgraf
\penalty-500\leavevmode\hss\endgraf
o49\ =\ 12\leavevmode\hss\endgraf
\endgroup
\penalty-1000
\par
\vskip 1 pt
\noindent
Since no eliminant\index{eliminant} given by a coordinate function satisfies
the hypotheses of the Shape Lemma, 
we took the eliminant with respect to the linear form
$y_{11} - y_{12} + y_{21} + y_{22}$.

This example is an instance of Lemma~3 of~\cite{SO:MPT00}.
These four points define a regular tetrahedron with volume
$V=9$ where each face has area $A=\sqrt{3^5}/2$ and each edge has length
$e=\sqrt{18}$.
That result guarantees that all 12 lines will be real when 
$e/2<r<A^2/3V$, which is the case above.

\section{Schubert calculus}
The classical Schubert calculus\index{Schubert calculus} of enumerative
geometry\index{enumerative geometry} concerns linear subspaces having
specified positions with respect to other, fixed subspaces. 
For instance, how many lines in ${\mathbb P}^3$ meet four given
lines? (See Example~\ref{ex:G22}.) 
More generally, let $1<r<n$ and suppose that we are given general linear
subspaces $L_1,\ldots,L_m$ of $k^n$ with $\dim L_i=n-r+1-l_i$.
When $l_1+\cdots+l_m=r(n-r)$, there will be a finite number
$d(r,n;l_1,\ldots,l_m)$ of $r$-planes in
$k^n$ which meet each $L_i$ non-trivially. 
This number may be computed using classical algorithms of Schubert and Pieri
(see~\cite{SO:MR48:2152}).

The condition on $r$-planes to meet a fixed $(n{-}r{+}1{-}l)$-plane
non-trivially is called a {\it (special) Schubert condition}, and we call 
the data $(r,n;l_1,\ldots,l_m)$ {\it (special) Schubert data}.
The {\it (special) Schubert calculus} concerns this class of enumerative
problems\index{Schubert calculus}\index{enumerative problem}.
We give two polynomial formulations of this special Schubert calculus,
consider their solutions over ${\mathbb R}$, and end with a question for
fields of arbitrary characteristic.

\subsection{Equations for the Grassmannian}\label{sec:grass}
The ambient space for the Schubert calculus is the
Grassmannian\index{Grassmannian}  
of $r$-planes in $k^n$, denoted ${\bf G}_{r,n}$.
For $H\in{\bf G}_{r,n}$, the $r$th exterior product of the embedding 
$H \rightarrow k^n$ gives a line
$$
  k\ \simeq\ \wedge^r H\ \longrightarrow\ \wedge^r k^{n}\ \simeq\ 
  k^{\binom{n}{r}}\,.
$$
This induces the Pl\"ucker embedding\index{Pl\"ucker embedding}  
${\bf G}_{r,n}\hookrightarrow{\mathbb P}^{\binom{n}{r}-1}$.
If $H$ is the row space of an $r$ by $n$ matrix, also written $H$, then the 
Pl\"ucker embedding sends $H$ to its vector of $\binom{n}{r}$ maximal
minors.
Thus the $r$-subsets of $\{0,\ldots,n{-}1\}$, 
${\mathbb Y}_{r,n}:={\tt subsets(n,r)}$,  
index Pl\"ucker coordinates\index{Pl\"ucker coordinate} of ${\bf G}_{r,n}$.
The Pl\"ucker ideal\index{Pl\"ucker ideal} of ${\bf G}_{r,n}$ is therefore the
ideal of algebraic relations among the maximal minors of a generic $r$ by $n$
matrix. 

We create the coordinate ring 
$k[p_\alpha\mid\alpha\in{\mathbb Y}_{2,5}]$ of ${\mathbb P}^9$ and the
Pl\"ucker ideal of ${\bf G}_{2,5}$.
The Grassmannian ${\bf G}_{r,n}$ of $r$-dimensional subspaces of $k^n$ is also
the Grassmannian of $r{-}1$-dimensional affine subspaces of 
${\mathbb P}^{n-1}$.
\Mtwo{}\index{\Mtwo{}} uses this alternative indexing scheme.
\par
\vskip 5 pt
\begingroup
\tteight
\baselineskip=10.01pt
\lineskip=0pt
\obeyspaces
i50\ :\ R\ =\ ZZ/101[apply(subsets(5,2),\ i\ ->\ p{\char`\_}i\ )];\leavevmode\hss\endgraf
\endgroup
\penalty-1000
\par
\vskip 1 pt
\noindent
\par
\vskip 5 pt
\begingroup
\tteight
\baselineskip=10.01pt
\lineskip=0pt
\obeyspaces
i51\ :\ I\ =\ Grassmannian(1,\ 4,\ R)\leavevmode\hss\endgraf
\penalty-500\leavevmode\hss\endgraf
o51\ =\ ideal\ (p\ \ \ \ \ \ p\ \ \ \ \ \ \ -\ p\ \ \ \ \ \ p\ \ \ \ \ \ \ +\ p\ \ \ \ \ \ p\ \ \ \ \ \ ,\ p\ \ \ \ \ \  $\cdot\cdot\cdot$\leavevmode\hss\endgraf
\ \ \ \ \ \ \ \ \ \ \ \ \ \ {\char`\{}2,\ 3{\char`\}}\ {\char`\{}1,\ 4{\char`\}}\ \ \ \ {\char`\{}1,\ 3{\char`\}}\ {\char`\{}2,\ 4{\char`\}}\ \ \ \ {\char`\{}1,\ 2{\char`\}}\ {\char`\{}3,\ 4{\char`\}}\ \ \ {\char`\{}2,\ 3{\char`\}} $\cdot\cdot\cdot$\leavevmode\hss\endgraf
\penalty-500\leavevmode\hss\endgraf
o51\ :\ Ideal\ of\ R\leavevmode\hss\endgraf
\endgroup
\penalty-1000
\par
\vskip 1 pt
\noindent
This projective variety has dimension 6 and degree 5
\par
\vskip 5 pt
\begingroup
\tteight
\baselineskip=10.01pt
\lineskip=0pt
\obeyspaces
i52\ :\ dim(Proj(R/I)),\ degree(I)\leavevmode\hss\endgraf
\penalty-500\leavevmode\hss\endgraf
o52\ =\ (6,\ 5)\leavevmode\hss\endgraf
\penalty-500\leavevmode\hss\endgraf
o52\ :\ Sequence\leavevmode\hss\endgraf
\endgroup
\penalty-1000
\par
\vskip 1 pt
\noindent

\begin{figure}[htb]
$$
 \epsfysize=2.16in \epsfbox{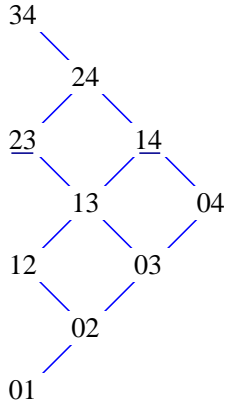}
$$\caption{${\mathbb Y}_{2,5}$\label{fig2}}
\end{figure}

This ideal has an important combinatorial 
structure~\cite[Example 11.9]{SO:Sturmfels_GBCP}. 
We write each $\alpha\in{\mathbb Y}_{r,n}$ as an increasing
sequence $\alpha\colon\alpha_1<\cdots<\alpha_r$. 
Given $\alpha,\beta\in{\mathbb Y}_{r,n}$, consider the 
two-rowed array with $\alpha$ written above $\beta$.
We say $\alpha\leq \beta$ if each column weakly increases.
If we sort the columns of an array with rows $\alpha$ and
$\beta$, then the first row is the {\it meet} $\alpha\wedge\beta$ 
(greatest lower bound) and the
second row the {\it join} $\alpha\vee\beta$ (least upper bound) of $\alpha$
and $\beta$. 
These definitions endow ${\mathbb Y}_{r,n}$ with the structure of a
distributive lattice.
Figure~\ref{fig2} shows ${\mathbb Y}_{2,5}$.

We give $k[p_\alpha]$ the degree reverse
lexicographic order, where we first order the variables $p_\alpha$ by
lexicographic order on their indices $\alpha$.

\begin{theorem}\label{PluckerIdeal}
The reduced Gr\"obner basis\index{Gr\"obner basis!reduced} of the Pl\"ucker
ideal with respect to this degree 
reverse lexicographic term order consists of quadratic
polynomials 
$$
  g(\alpha,\beta)\quad=\quad
  p_\alpha\cdot p_\beta \ -\  p_{\alpha\vee\beta}\cdot p_{\alpha\wedge\beta} 
  \ +\ \hbox{lower terms in $\prec$}\,,
$$
for each incomparable pair $\alpha,\beta$ in
${\mathbb Y}_{r,n}$,
and all lower terms $\lambda p_\gamma\cdot p_\delta$ in $g(\alpha,\beta)$
satisfy $\gamma\leq \alpha\wedge\beta$ and $\alpha\vee\beta\leq \delta$.
\end{theorem}

The form of this Gr\"obner basis implies that the standard monomials are 
the sortable monomials, those $p_\alpha p_\beta\cdots p_\gamma$ with 
$\alpha\leq\beta\leq\cdots\leq\gamma$.
Thus the Hilbert function\index{Hilbert function} of ${\bf G}_{r,n}$ may be
expressed in terms of the combinatorics of ${\mathbb Y}_{r,n}$.
For instance, the dimension of ${\bf G}_{r,n}$ is the rank of 
${\mathbb Y}_{r,n}$, and its degree is the number of maximal chains.
From Figure~\ref{fig2}, these are 6 and 5 for ${\mathbb Y}_{2,5}$,
confirming our previous calculations.

Since the generators $g(\alpha,\beta)$ are linearly independent, this
Gr\"obner basis is also a minimal generating set for the ideal.
The displayed generator in {\tt o51}, 
$$
  p_{\{2,3\}}p_{\{1,4\}}\ -\ p_{\{1,3\}}p_{\{2,4\}}\ -\
  p_{\{1,2\}}p_{\{3,4\}}\ ,
$$
is $g(23, 14)$, and corresponds to the underlined incomparable pair in
Figure~\ref{fig2}. 
Since there are 5 such incomparable pairs, the Gr\"obner basis has 5
generators.
As ${\bf G}_{2,5}$ has codimension 3, it is not a complete
intersection\index{Grassmannian!not a complete intersection}. 
This shows how the general enumerative problem from the Schubert calculus
gives rise to an overdetermined system of equations\index{polynomial
equations!overdetermined} in this global 
formulation. 
\medskip

The Grassmannian\index{Grassmannian!local coordinates} has a useful system of
local coordinates given by ${\rm Mat}_{r,n-r}$ as follows
\begin{equation}\label{eq:local}
  Y\ \in {\rm Mat}_{r,n-r}\ \longmapsto\ 
  {\rm row space}\ [ I_r : Y ]\ \in\ {\bf G}_{r,n}\,.
\end{equation}

Let $L$ be a ($n-r+1-l$)-plane in $k^n$ which is the row space of
a $n-r+1-l$ by $n$ matrix, also written $L$.
Then $L$ meets $X\in{\bf G}_{r,n}$ non-trivially if
$$
  \mbox{maximal minors of }\ 
  \left[\begin{array}{c}L\\X\end{array}\right]\ =\ 0\,.
$$
Laplace expansion of each minor along the rows of $X$ gives a linear
equation in the Pl\"ucker coordinates.
In the local coordinates (substituting $[I_r:Y]$ for $X$), we obtain
multilinear equations of degree $\min\{r,n-r\}$.
These equations generate a prime ideal of codimension $l$.

Suppose each $l_i=1$ in our enumerative problem.
Then in the Pl\"ucker coordinates, we have the Pl\"ucker ideal of 
${\bf G}_{r,n}$
together with $r(n-r)$ linear equations, one for each
$(n{-}r)$-plane $L_i$.
By Theorem~\ref{PluckerIdeal}, the Pl\"ucker ideal has a square-free initial
ideal\index{initial ideal!square-free}, and so the Gr\"obner
homotopy\index{homotopy!Gr\"obner} of 
Example~\ref{example:Groebner} may be 
used to solve this enumerative problem.

\begin{example}\label{ex:G22}
${\bf G}_{2,4}\subset{\mathbb P}^5$ has equation
\begin{equation}\label{eq:G22}
   p_{\{1,2\}}p_{\{0,3\}}-p_{\{1,3\}}p_{\{0,2\}}+ p_{\{2,3\}}p_{\{0,1\}}
   \ =\ 0\,.
\end{equation}
The condition for $H\in{\bf G}_{2,4}$ to meet a 2-plane $L$ is the
vanishing of 
\begin{equation}\label{eq:hypersurface}
    p_{\{1,2\}}L_{34}-p_{\{1,3\}}L_{24}+p_{\{2,3\}}L_{14} 
  + p_{\{1,4\}}L_{23}-p_{\{2,4\}}L_{13}+p_{\{3,4\}}L_{12}\,,
\end{equation}
where $L_{ij}$ is the $(i,j)$th maximal minor of $L$.

If $l_1=\cdots=l_4=1$, we have 5 equations in ${\mathbb P}^5$, one quadratic
and 4 linear, and so by B\'ezout's Theorem\index{B\'ezout Theorem} there are
two 2-planes in $k^4$ that meet 4 general 2-planes non-trivially. 
This means that there are 2 lines in ${\mathbb P}^3$ meeting 4 general lines. 
In local coordinates, (\ref{eq:hypersurface}) becomes
$$
    L_{34}-L_{14}y_{11}+L_{13}y_{12}-L_{24}y_{21} 
  + L_{23}y_{22} + L_{12}(y_{11}y_{22}-y_{12}y_{21})\,.
$$
This polynomial has the form of the last specialization in
Example~\ref{ex:one}. 
\end{example}

\subsection{Reality in the Schubert calculus}\label{sec:shapiro}
Like the other enumerative problems we have discussed, enumerative problems
in the special Schubert calculus\index{Schubert calculus} are fully
real\index{enumerative problem!fully real} in that all solutions can be 
real~\cite{SO:So99a}. 
That is, given any Schubert data $(r,n;l_1,\ldots,l_m)$, there exist 
subspaces $L_1,\ldots,L_m\subset{\mathbb R}^n$ such that each of the
$d(r,n;l_1,\ldots,l_m)$ $r$-planes that meet each $L_i$ are themselves
real.

This result gives some idea of which choices of the $L_i$ give
all $r$-planes real.
Let $\gamma$ be a fixed rational normal curve in ${\mathbb R}^n$.
Then the $L_i$ are linear subspaces osculating $\gamma$.
More concretely, suppose that $\gamma$ is the standard rational normal
curve\index{rational normal curve},
$\gamma(s) = (1, s, s^2, \ldots, s^{n-1})$.
Then the $i$-plane 
$L_i(s):=\langle \gamma(s),\gamma'(s),\ldots,\gamma^{(i-1)}(s)\rangle$ 
osculating $\gamma$ at $\gamma(s)$ is the row space
of the matrix given by {\tt oscPlane(i, n, s)}.
\par
\vskip 5 pt
\begingroup
\tteight
\baselineskip=10.01pt
\lineskip=0pt
\obeyspaces
i53\ :\ oscPlane\ =\ (i,\ n,\ s)\ ->\ (\leavevmode\hss\endgraf
\ \ \ \ \ \ \ \ \ \ \ gamma\ :=\ matrix\ {\char`\{}toList\ apply(1..n,\ i\ ->\ s{\char`\^}(i-1)){\char`\}};\leavevmode\hss\endgraf
\ \ \ \ \ \ \ \ \ \ \ L\ :=\ gamma;\leavevmode\hss\endgraf
\ \ \ \ \ \ \ \ \ \ \ j\ :=\ 0;\leavevmode\hss\endgraf
\ \ \ \ \ \ \ \ \ \ \ while\ j\ <\ i-1\ do\ (gamma\ =\ diff(s,\ gamma);\ \leavevmode\hss\endgraf
\ \ \ \ \ \ \ \ \ \ \ \ \ \ \ \ L\ =\ L\ ||\ gamma;\leavevmode\hss\endgraf
\ \ \ \ \ \ \ \ \ \ \ \ \ \ \ \ j\ =\ j+1);\leavevmode\hss\endgraf
\ \ \ \ \ \ \ \ \ \ \ L);\leavevmode\hss\endgraf
\endgroup
\penalty-1000
\par
\vskip 1 pt
\noindent
\par
\vskip 5 pt
\begingroup
\tteight
\baselineskip=10.01pt
\lineskip=0pt
\obeyspaces
i54\ :\ QQ[s];\ oscPlane(3,\ 6,\ s)\leavevmode\hss\endgraf
\penalty-500\leavevmode\hss\endgraf
o55\ =\ |\ 1\ s\ s2\ s3\ \ s4\ \ \ s5\ \ \ |\leavevmode\hss\endgraf
\ \ \ \ \ \ |\ 0\ 1\ 2s\ 3s2\ 4s3\ \ 5s4\ \ |\leavevmode\hss\endgraf
\ \ \ \ \ \ |\ 0\ 0\ 2\ \ 6s\ \ 12s2\ 20s3\ |\leavevmode\hss\endgraf
\penalty-500\leavevmode\hss\endgraf
\ \ \ \ \ \ \ \ \ \ \ \ \ \ \ \ \ \ \ 3\ \ \ \ \ \ \ \ \ \ \ \ 6\leavevmode\hss\endgraf
o55\ :\ Matrix\ QQ\ [s]\ \ <---\ QQ\ [s]\leavevmode\hss\endgraf
\endgroup
\penalty-1000
\par
\vskip 1 pt
\noindent
(In {\tt o55}, the exponents of $s$ are displayed in line: $s^2$ is written
{\tt s2}.
\Mtwo{}\index{\Mtwo{}} uses this notational convention to display
matrices efficiently.)

\begin{theorem}[\cite{SO:So99a}]\label{thm:special-reality}
For any Schubert data $(r,n;l_1,\ldots,l_m)$, {\bf there exist} real numbers
$s_1,s_2,\ldots,s_m$ such that there are $d(r,n;l_1,\ldots,l_m)$
$r$-planes that meet each osculating plane $L_i(s_i)$, and all are real.
\end{theorem}

The inspiration for looking at subspaces osculating the rational normal
curve\index{rational normal curve} to
study real enumerative geometry\index{enumerative geometry!real} for the
Schubert calculus\index{Schubert calculus} is the following very interesting 
conjecture of Boris Shapiro\index{Shapiro, B.} and Michael
Shapiro\index{Shapiro, M.}, or more accurately,
extensive computer experimentation based upon their
conjecture~\cite{SO:RS98,SO:So_shap-www,SO:So00b,SO:Ver00}.
\medskip

\noindent{\bf Shapiros's Conjecture\index{Shapiros's Conjecture}. }
{\it 
For any Schubert data $(r,n;l_1,\ldots,l_m)$ and {\bf for all} real numbers
$s_1,s_2,\ldots,s_m$ there are $d(r,n;l_1,\ldots,l_m)$
$r$-planes that meet each osculating plane $L_i(s_i)$, and all are real.
}\medskip

In addition to Theorem~\ref{thm:special-reality}, (which replaces the 
quantifier {\it for all}\/ by  {\it there exist}), the strongest evidence for
this Conjecture is the following result of Eremenko\index{Eremenko, A.} and
Gabrielov\index{Gabrielov, A.}~\cite{SO:EG00}. 

\begin{theorem}
Shapiros's Conjecture is true when either $r$ or $n-r$ is $2$.
\end{theorem}

We test an example of this conjecture for the Schubert data
$(3,6;1^3,2^3)$, (where $a^b$ is $a$ repeated $b$ times).
The algorithms of the Schubert calculus predict that $d(3,6;1^3,2^3)=6$.
The function {\tt spSchub(r, L, P)} computes the ideal of $r$-planes meeting
the row space of $L$ in the Pl\"ucker coordinates $P_\alpha$.
\par
\vskip 5 pt
\begingroup
\tteight
\baselineskip=10.01pt
\lineskip=0pt
\obeyspaces
i56\ :\ spSchub\ =\ (r,\ L,\ P)\ ->\ (\leavevmode\hss\endgraf
\ \ \ \ \ \ \ \ \ \ \ I\ :=\ ideal\ apply(subsets(numgens\ source\ L,\ \leavevmode\hss\endgraf
\ \ \ \ \ \ \ \ \ \ \ \ \ \ \ \ \ \ \ \ \ \ \ \ \ \ \ \ r\ +\ numgens\ target\ L),\ S\ ->\ \leavevmode\hss\endgraf
\ \ \ \ \ \ \ \ \ \ \ \ \ \ \ \ fold((sum,\ U)\ ->\ sum\ +\leavevmode\hss\endgraf
\ \ \ \ \ \ \ \ \ \ \ \ \ \ \ \ \ fold((term,i)\ ->\ term*(-1){\char`\^}i,\ P{\char`\_}(S{\char`\_}U)\ *\ det(\leavevmode\hss\endgraf
\ \ \ \ \ \ \ \ \ \ \ \ \ \ \ \ \ \ submatrix(L,\ sort\ toList(set(S)\ -\ set(S{\char`\_}U)))),\ U),\ \leavevmode\hss\endgraf
\ \ \ \ \ \ \ \ \ \ \ \ \ \ \ \ \ \ \ \ \ 0,\ subsets({\char`\#}S,\ r))));\leavevmode\hss\endgraf
\endgroup
\penalty-1000
\par
\vskip 1 pt
\noindent
We are working in the Grassmannian of 3-planes in 
${\mathbb C}^6$.
\par
\vskip 5 pt
\begingroup
\tteight
\baselineskip=10.01pt
\lineskip=0pt
\obeyspaces
i57\ :\ R\ =\ QQ[apply(subsets(6,3),\ i\ ->\ p{\char`\_}i\ )];\leavevmode\hss\endgraf
\endgroup
\penalty-1000
\par
\vskip 1 pt
\noindent
The ideal $I$ consists of the
special Schubert conditions for the 3-planes to meet the 3-planes osculating
the rational normal curve at the points 1, 2, and 3, and to also meet the
2-planes osculating at 4, 5, and 6,
together with the Pl\"ucker ideal {\tt Grassmannian(2, 5, R)}.
Since this is a 1-dimensional homogeneous ideal, we add the linear form 
{\tt p{\char`\_}\{0,1,5\} - 1} to make the ideal
zero-dimensional\index{ideal!zero-dimensional}.
As before, {\tt Grassmannian(2, 5, R)} creates the Pl\"ucker ideal of 
${\bf G}_{3,6}$.
\par
\vskip 5 pt
\begingroup
\tteight
\baselineskip=10.01pt
\lineskip=0pt
\obeyspaces
i58\ :\ I\ =\ fold((J,\ i)\ ->\ J\ +\leavevmode\hss\endgraf
\ \ \ \ \ \ \ \ \ \ \ \ spSchub(3,\ substitute(oscPlane(3,\ 6,\ s),\ {\char`\{}s=>\ 1+i{\char`\}}),\ p)\ +\leavevmode\hss\endgraf
\ \ \ \ \ \ \ \ \ \ \ \ spSchub(3,\ substitute(oscPlane(2,\ 6,\ s),\ {\char`\{}s=>\ 4+i{\char`\}}),\ p),\ \leavevmode\hss\endgraf
\ \ \ \ \ \ \ \ \ \ \ \ Grassmannian(2,\ 5,\ R),\ {\char`\{}0,1,2{\char`\}})\ +\ \leavevmode\hss\endgraf
\ \ \ \ \ \ \ \ \ \ \ ideal\ (p{\char`\_}{\char`\{}0,1,5{\char`\}}\ -\ 1);\leavevmode\hss\endgraf
\penalty-500\leavevmode\hss\endgraf
o58\ :\ Ideal\ of\ R\leavevmode\hss\endgraf
\endgroup
\penalty-1000
\par
\vskip 1 pt
\noindent
This has dimension 0 and degree 6, in agreement with the Schubert calculus.
\par
\vskip 5 pt
\begingroup
\tteight
\baselineskip=10.01pt
\lineskip=0pt
\obeyspaces
i59\ :\ dim\ I,\ degree\ I\leavevmode\hss\endgraf
\penalty-500\leavevmode\hss\endgraf
o59\ =\ (0,\ 6)\leavevmode\hss\endgraf
\penalty-500\leavevmode\hss\endgraf
o59\ :\ Sequence\leavevmode\hss\endgraf
\endgroup
\penalty-1000
\par
\vskip 1 pt
\noindent
As expected, all roots are real.
\par
\vskip 5 pt
\begingroup
\tteight
\baselineskip=10.01pt
\lineskip=0pt
\obeyspaces
i60\ :\ A\ =\ R/I;\ numRealSturm(eliminant(p{\char`\_}{\char`\{}2,3,4{\char`\}},\ QQ[Z]))\leavevmode\hss\endgraf
\penalty-500\leavevmode\hss\endgraf
o61\ =\ 6\leavevmode\hss\endgraf
\endgroup
\penalty-1000
\par
\vskip 1 pt
\noindent
There have been many checked instances of this
conjecture~\cite{SO:So_shap-www,SO:So00b,SO:Ver00}, and it has some
geometrically interesting generalizations~\cite{SO:So_flags}.

The question remains for which numbers $0\leq d\leq d(r,n;l_1,\ldots,l_m)$ do
there exist real planes $L_i$ with $d(r,n;l_1,\ldots,l_m)$
$r$-planes meeting each $L_i$, and exactly $d$ of them are real.
Besides Theorem~\ref{thm:special-reality} and the obvious parity condition,
nothing is known in general.
In every known case, every possibility occurs---which is not the case in all
enumerative problems, even those that are fully real\index{enumerative
problem!fully real}\footnote{For example, of
the 12 rational plane cubics containing 8 real points in ${\mathbb P}^2$,
either 8, 10 or 12 can be real, and there are 8 points with all 12
real~\cite[Proposition 4.7.3]{SO:DeKh00}.}.
Settling this (for $d=0$) has implications for linear systems
theory~\cite{SO:RS98}. 

\subsection{Transversality in the Schubert calculus}
A basic principle of the classical Schubert calculus\index{Schubert calculus}
is that the intersection 
number $d(r,n;l_1,\ldots,l_m)$ has enumerative significance---that is, for
general linear subspaces $L_i$, all solutions appear with multiplicity 1.
This basic principle is not known to hold in general.
For fields of characteristic zero, Kleiman's Transversality
Theorem~\cite{SO:MR50:13063} establishes this principle.
When $r$ or $n{-}r$ is 2, then Theorem~E of~\cite{SO:So97a} establishes this
principle in arbitrary characteristic.
We conjecture that this principle holds in general; that is, for arbitrary
infinite fields and any Schubert data, if the planes $L_i$ are in general
position, then the resulting zero-dimensional ideal is
radical\index{ideal!radical}. 

We test this conjecture on the enumerative problem of 
Section~\ref{sec:shapiro}, which is not covered by 
Theorem~E of~\cite{SO:So97a}.
The function {\tt testTransverse(F)} tests transversality
for this enumerative problem, for a given field $F$. 
It does this by first computing the ideal of the enumerative problem using 
random planes $L_i$.
\par
\vskip 5 pt
\begingroup
\tteight
\baselineskip=10.01pt
\lineskip=0pt
\obeyspaces
i62\ :\ randL\ =\ (R,\ n,\ r,\ l)\ ->\ \leavevmode\hss\endgraf
\ \ \ \ \ \ \ \ \ \ \ \ \ \ \ \ matrix\ table(n-r+1-l,\ n,\ (i,\ j)\ ->\ random(0,\ R));\leavevmode\hss\endgraf
\endgroup
\penalty-1000
\par
\vskip 1 pt
\noindent
and the Pl\"ucker ideal of the Grassmannian ${\bf G}_{3,6}$
 {\tt Grassmannian(2, 5, R)}.)
Then it adds a random (inhomogeneous) linear relation 
 {\tt 1 + random(1, R)} to make the ideal zero-dimensional for generic $L_i$. 
When this ideal is zero dimensional and has degree 6 (the expected degree), it
computes the characteristic polynomial {\tt g} of a generic linear form.
If {\tt g} has no multiple roots, {\tt 1 == gcd(g, diff(Z, g))}, 
then the Shape Lemma\index{Shape Lemma}
guarantees that the ideal was radical.
{\tt testTransverse} exits either when it computes a radical ideal,
or after {\tt limit} iterations (which is set to 5 for these examples), and
prints the return status. 
\par
\vskip 5 pt
\begingroup
\tteight
\baselineskip=10.01pt
\lineskip=0pt
\obeyspaces
i63\ :\ testTransverse\ =\ F\ ->\ (\leavevmode\hss\endgraf
\ \ \ \ \ \ \ \ \ \ \ \ R\ :=\ F[apply(subsets(6,\ 3),\ i\ ->\ q{\char`\_}i\ )];\leavevmode\hss\endgraf
\ \ \ \ \ \ \ \ \ \ \ \ continue\ :=\ true;\leavevmode\hss\endgraf
\ \ \ \ \ \ \ \ \ \ \ \ j\ :=\ 0;\ \ \leavevmode\hss\endgraf
\ \ \ \ \ \ \ \ \ \ \ \ limit\ :=\ 5;\leavevmode\hss\endgraf
\ \ \ \ \ \ \ \ \ \ \ \ while\ continue\ and\ (j\ <\ limit)\ do\ (\leavevmode\hss\endgraf
\ \ \ \ \ \ \ \ \ \ \ \ \ \ \ \ \ j\ =\ j\ +\ 1;\leavevmode\hss\endgraf
\ \ \ \ \ \ \ \ \ \ \ \ \ \ \ \ \ I\ :=\ fold((J,\ i)\ ->\ J\ +\ \leavevmode\hss\endgraf
\ \ \ \ \ \ \ \ \ \ \ \ \ \ \ \ \ \ \ \ \ \ \ \ \ \ \ spSchub(3,\ randL(R,\ 6,\ 3,\ 1),\ q)\ +\leavevmode\hss\endgraf
\ \ \ \ \ \ \ \ \ \ \ \ \ \ \ \ \ \ \ \ \ \ \ \ \ \ \ spSchub(3,\ randL(R,\ 6,\ 3,\ 2),\ q),\leavevmode\hss\endgraf
\ \ \ \ \ \ \ \ \ \ \ \ \ \ \ \ \ \ \ \ \ \ \ \ \ \ \ Grassmannian(2,\ 5,\ R)\ +\ \leavevmode\hss\endgraf
\ \ \ \ \ \ \ \ \ \ \ \ \ \ \ \ \ \ \ \ \ \ \ \ \ \ \ ideal\ (1\ +\ random(1,\ R)),\leavevmode\hss\endgraf
\ \ \ \ \ \ \ \ \ \ \ \ \ \ \ \ \ \ \ \ \ \ \ \ \ \ \ {\char`\{}0,\ 1,\ 2{\char`\}});\leavevmode\hss\endgraf
\ \ \ \ \ \ \ \ \ \ \ \ \ \ \ \ \ if\ (dim\ I\ ==\ 0)\ and\ (degree\ I\ ==\ 6)\ then\ (\leavevmode\hss\endgraf
\ \ \ \ \ \ \ \ \ \ \ \ \ \ \ \ \ lin\ :=\ promote(random(1,\ R),\ (R/I));\leavevmode\hss\endgraf
\ \ \ \ \ \ \ \ \ \ \ \ \ \ \ \ \ g\ :=\ charPoly(lin,\ Z);\leavevmode\hss\endgraf
\ \ \ \ \ \ \ \ \ \ \ \ \ \ \ \ \ continue\ =\ not(1\ ==\ gcd(g,\ diff(Z,\ g)));\leavevmode\hss\endgraf
\ \ \ \ \ \ \ \ \ \ \ \ \ \ \ \ \ ));\leavevmode\hss\endgraf
\ \ \ \ \ \ \ \ \ \ \ \ if\ continue\ then\ <<\ "Failed\ for\ the\ prime\ "\ <<\ char\ F\ <<\ \leavevmode\hss\endgraf
\ \ \ \ \ \ \ \ \ \ \ \ \ \ \ "\ with\ "\ <<\ j\ <<\ "\ iterations"\ <<\ endl;\leavevmode\hss\endgraf
\ \ \ \ \ \ \ \ \ \ \ \ if\ not\ continue\ then\ <<\ "Succeeded\ for\ the\ prime\ "\ <<\leavevmode\hss\endgraf
\ \ \ \ \ \ \ \ \ \ \ \ \ \ \ \ char\ F\ <<\ "\ in\ "\ <<\ j\ <<\ "\ iteration(s)"\ <<\ endl;\leavevmode\hss\endgraf
\ \ \ \ \ \ \ \ \ \ \ \ );\leavevmode\hss\endgraf
\endgroup
\penalty-1000
\par
\vskip 1 pt
\noindent
Since 5 iterations do not show transversality for ${\mathbb F}_2$,
\par
\vskip 5 pt
\begingroup
\tteight
\baselineskip=10.01pt
\lineskip=0pt
\obeyspaces
i64\ :\ testTransverse(ZZ/2);\leavevmode\hss\endgraf
Failed\ for\ the\ prime\ 2\ with\ 5\ iterations\leavevmode\hss\endgraf
\endgroup
\penalty-1000
\par
\vskip 1 pt
\noindent
we can test transversality in characteristic 2 using the field with
four elements, ${\mathbb F}_4=$ {\tt GF 4}.
\par
\vskip 5 pt
\begingroup
\tteight
\baselineskip=10.01pt
\lineskip=0pt
\obeyspaces
i65\ :\ testTransverse(GF\ 4);\leavevmode\hss\endgraf
Succeeded\ for\ the\ prime\ 2\ in\ 3\ iteration(s)\leavevmode\hss\endgraf
\endgroup
\penalty-1000
\par
\vskip 1 pt
\noindent
We do find transversality for ${\mathbb F}_7$.
\par
\vskip 5 pt
\begingroup
\tteight
\baselineskip=10.01pt
\lineskip=0pt
\obeyspaces
i66\ :\ testTransverse(ZZ/7);\leavevmode\hss\endgraf
Succeeded\ for\ the\ prime\ 7\ in\ 2\ iteration(s)\leavevmode\hss\endgraf
\endgroup
\penalty-1000
\par
\vskip 1 pt
\noindent

We have tested transversality for all primes less than 100 in every
enumerative problem involving Schubert conditions on 3-planes in $k^6$.
These include the problem above as well as the problem of 42 3-planes meeting
9 general 3-planes.\footnote{After this was written, we discovered an
elementary proof of transversality for the enumerative problems 
$d(r,n;1^{r(n-r)})$, where the conditions are all
codimension~1~\cite{SO:So_trans}.}

\section{The 12 lines: reprise}
The enumerative problems of Section~\ref{sec:enumerative} were formulated
in local coordinates~(\ref{eq:local}) for the Grassmannian of lines in
${\mathbb P}^3$ (Grassmannian of 2-dimensional subspaces in $k^4$).
When we formulate the problem of Section~\ref{sec:12lines} in the global
Pl\"ucker coordinates\index{Pl\"ucker coordinate} of Section~\ref{sec:grass},
we find some interesting phenomena.
We also consider some related enumerative problems\index{enumerative problem}.

\subsection{Global formulation}\label{sec:global}
A quadratic form\index{quadratic form} $q$ on a vector space $V$ over a field
$k$ not of characteristic 2 is given by
$q({\bf u})=(\varphi({\bf u}),{\bf u})$, where $\varphi\colon V\to V^*$ is
a {\it symmetric} linear map, that is 
$(\varphi({\bf u}),{\bf v})=(\varphi({\bf v}),{\bf u})$.
Here, $V^*$ is the linear dual of $V$ and $(\,\cdot\;,\,\cdot\,)$ is the
pairing $V\otimes V^*\to k$. 
The map $\varphi$ induces a quadratic form $\wedge^rq$ on the $r$th exterior
power $\wedge^rV$ of $V$ through the symmetric map 
$\wedge^r\varphi\colon \wedge^rV\to\wedge^rV^*=(\wedge^rV)^*$.
The action of $\wedge^rV^*$ on $\wedge^rV$ is given by
\begin{equation}\label{eq:wedge}
  ({\bf x}_1\wedge{\bf x}_2\wedge\cdots\wedge{\bf x}_r,\ 
   {\bf y\!}_1\wedge{\bf y}\!_2\wedge\cdots\wedge{\bf y}\!_r)\ =\ 
   \det|({\bf x}_i, {\bf y}\!_j)|\,,
\end{equation}
where ${\bf x}_i\in V^*$ and ${\bf y}\!_j\in V$.

When we fix isomorphisms $V\simeq k^n\simeq V^*$, the map $\varphi$ is given
by a symmetric $n\times n$ matrix $M$ as in Section~\ref{sec:12lines}.
Suppose $r=2$.
Then for ${\bf u},{\bf v}\in k^n$,
$$
  \wedge^2q({\bf u}\wedge{\bf v})\ = 
              \det\left[\begin{array}{cc}
                    {\bf u}M{\bf u}^t &{\bf u}M{\bf v}^t\\
                    {\bf v}M{\bf u}^t &{\bf v}M{\bf v}^t
                  \end{array}\right]\ ,
$$
which is Equation~(\ref{eq:sphere}) of  Section~\ref{sec:12lines}.

\begin{proposition}\label{prop:tangent_line}
  A line $\ell$ is tangent to a quadric ${\mathcal V}(q)$ in 
  ${\mathbb P}^{n-1}$ if and only if its Pl\"ucker 
  coordinate\index{Pl\"ucker coordinate}  
  $\wedge^2\ell\in{\mathbb P}^{\binom{n}{2}-1}$ lies on the quadric
  ${\mathcal V}(\wedge^2q)$.
\end{proposition}
       
Thus the Pl\"ucker coordinates for the set of lines tangent to 4 general
quadrics in ${\mathbb P}^3$ satisfy 5 quadratic equations:
The single Pl\"ucker relation~(\ref{eq:G22})
together with one quadratic equation for each quadric.
Thus we expect the B\'ezout number\index{B\'ezout number} of $2^5=32$ such
lines. 
We check this.

The procedure {\tt randomSymmetricMatrix(R, n)}
generates a random symmetric $n\times n$ matrix with entries in 
the base ring of $R$.
\par
\vskip 5 pt
\begingroup
\tteight
\baselineskip=10.01pt
\lineskip=0pt
\obeyspaces
i67\ :\ randomSymmetricMatrix\ =\ (R,\ n)\ ->\ (\leavevmode\hss\endgraf
\ \ \ \ \ \ \ \ \ \ entries\ :=\ new\ MutableHashTable;\leavevmode\hss\endgraf
\ \ \ \ \ \ \ \ \ \ scan(0..n-1,\ i\ ->\ scan(i..n-1,\ j\ ->\ \leavevmode\hss\endgraf
\ \ \ \ \ \ \ \ \ \ \ \ \ \ \ \ \ \ \ \ \ \ \ entries{\char`\#}(i,\ j)\ =\ random(0,\ R)));\leavevmode\hss\endgraf
\ \ \ \ \ \ \ \ \ \ matrix\ table(n,\ n,\ (i,\ j)\ ->\ if\ i\ >\ j\ then\ \leavevmode\hss\endgraf
\ \ \ \ \ \ \ \ \ \ \ \ \ \ \ \ \ \ \ \ \ \ \ entries{\char`\#}(j,\ i)\ else\ entries{\char`\#}(i,\ j))\leavevmode\hss\endgraf
\ \ \ \ \ \ \ \ \ \ );\leavevmode\hss\endgraf
\endgroup
\penalty-1000
\par
\vskip 1 pt
\noindent
The procedure {\tt tangentEquation(r, R, M)} gives the equation in Pl\"ucker
coordinates for a point in ${\mathbb P}^{\binom{n}{r}-1}$ to be
isotropic with respect to the bilinear form\index{bilinear form}  $\wedge^rM$
({\tt R} is assumed to be the coordinate ring of 
${\mathbb P}^{\binom{n}{r}-1}$).
This is the equation for an $r$-plane to be tangent to the quadric
associated to $M$.
\par
\vskip 5 pt
\begingroup
\tteight
\baselineskip=10.01pt
\lineskip=0pt
\obeyspaces
i68\ :\ tangentEquation\ =\ (r,\ R,\ M)\ ->\ (\leavevmode\hss\endgraf
\ \ \ \ \ \ \ \ \ \ \ g\ :=\ matrix\ {\char`\{}gens(R){\char`\}};\leavevmode\hss\endgraf
\ \ \ \ \ \ \ \ \ \ \ (entries(g\ *\ exteriorPower(r,\ M)\ *\ transpose\ g)){\char`\_}0{\char`\_}0\leavevmode\hss\endgraf
\ \ \ \ \ \ \ \ \ \ \ );\leavevmode\hss\endgraf
\endgroup
\penalty-1000
\par
\vskip 1 pt
\noindent
We construct the ideal of lines tangent to 4 general quadrics in
${\mathbb P}^3$.
\par
\vskip 5 pt
\begingroup
\tteight
\baselineskip=10.01pt
\lineskip=0pt
\obeyspaces
i69\ :\ R\ =\ QQ[apply(subsets(4,\ 2),\ i\ ->\ p{\char`\_}i\ )];\leavevmode\hss\endgraf
\endgroup
\penalty-1000
\par
\vskip 1 pt
\noindent
\par
\vskip 5 pt
\begingroup
\tteight
\baselineskip=10.01pt
\lineskip=0pt
\obeyspaces
i70\ :\ I\ =\ Grassmannian(1,\ 3,\ R)\ +\ ideal\ apply(0..3,\ i\ ->\ \leavevmode\hss\endgraf
\ \ \ \ \ \ \ \ \ \ \ tangentEquation(2,\ R,\ randomSymmetricMatrix(R,\ 4)));\leavevmode\hss\endgraf
\penalty-500\leavevmode\hss\endgraf
o70\ :\ Ideal\ of\ R\leavevmode\hss\endgraf
\endgroup
\penalty-1000
\par
\vskip 1 pt
\noindent
As expected, this ideal has dimension 0 and degree 32.
\par
\vskip 5 pt
\begingroup
\tteight
\baselineskip=10.01pt
\lineskip=0pt
\obeyspaces
i71\ :\ dim\ Proj(R/I),\ degree\ I\leavevmode\hss\endgraf
\penalty-500\leavevmode\hss\endgraf
o71\ =\ (0,\ 32)\leavevmode\hss\endgraf
\penalty-500\leavevmode\hss\endgraf
o71\ :\ Sequence\leavevmode\hss\endgraf
\endgroup
\penalty-1000
\par
\vskip 1 pt
\noindent

\subsection{Lines tangent to 4 spheres}\label{sec:tangent_lines}
That calculation raises the following question:
In Section~\ref{sec:12lines}, why did we obtain only 12 lines tangent to 4
spheres\index{sphere}? 
To investigate this, we generate the global ideal of lines tangent to the 
spheres of Section~\ref{sec:12lines}.
\par
\vskip 5 pt
\begingroup
\tteight
\baselineskip=10.01pt
\lineskip=0pt
\obeyspaces
i72\ :\ I\ =\ Grassmannian(1,\ 3,\ R)\ +\ \leavevmode\hss\endgraf
\ \ \ \ \ \ \ \ \ \ \ \ \ \ ideal\ (tangentEquation(2,\ R,\ Sphere(0,0,0,5)),\leavevmode\hss\endgraf
\ \ \ \ \ \ \ \ \ \ \ \ \ \ \ \ \ \ \ \ \ tangentEquation(2,\ R,\ Sphere(4,1,1,5)),\leavevmode\hss\endgraf
\ \ \ \ \ \ \ \ \ \ \ \ \ \ \ \ \ \ \ \ \ tangentEquation(2,\ R,\ Sphere(1,4,1,5)),\leavevmode\hss\endgraf
\ \ \ \ \ \ \ \ \ \ \ \ \ \ \ \ \ \ \ \ \ tangentEquation(2,\ R,\ Sphere(1,1,4,5)));\leavevmode\hss\endgraf
\penalty-500\leavevmode\hss\endgraf
o72\ :\ Ideal\ of\ R\leavevmode\hss\endgraf
\endgroup
\penalty-1000
\par
\vskip 1 pt
\noindent
We compute the dimension and degree of ${\mathcal V}(I)$.
\par
\vskip 5 pt
\begingroup
\tteight
\baselineskip=10.01pt
\lineskip=0pt
\obeyspaces
i73\ :\ dim\ Proj(R/I),\ degree\ I\leavevmode\hss\endgraf
\penalty-500\leavevmode\hss\endgraf
o73\ =\ (1,\ 4)\leavevmode\hss\endgraf
\penalty-500\leavevmode\hss\endgraf
o73\ :\ Sequence\leavevmode\hss\endgraf
\endgroup
\penalty-1000
\par
\vskip 1 pt
\noindent
The ideal is not zero dimensional\index{ideal!zero-dimensional}; there is an
extraneous one-dimensional 
component of zeroes with degree 4.
Since we found 12 lines in
Section~\ref{sec:12lines} using the local coordinates~(\ref{eq:local}),
the extraneous component must lie in the complement of that coordinate patch,
which is defined by the vanishing of the first Pl\"ucker coordinate, $p_{\{0,1\}}$.
We saturate\index{saturate} $I$ by $p_{\{0,1\}}$ to obtain the desired lines.
\par
\vskip 5 pt
\begingroup
\tteight
\baselineskip=10.01pt
\lineskip=0pt
\obeyspaces
i74\ :\ Lines\ =\ saturate(I,\ ideal\ (p{\char`\_}{\char`\{}0,1{\char`\}}));\leavevmode\hss\endgraf
\penalty-500\leavevmode\hss\endgraf
o74\ :\ Ideal\ of\ R\leavevmode\hss\endgraf
\endgroup
\penalty-1000
\par
\vskip 1 pt
\noindent
This ideal does have dimension 0 and degree 12, so we have recovered the
zeroes of Section~\ref{sec:12lines}.
\par
\vskip 5 pt
\begingroup
\tteight
\baselineskip=10.01pt
\lineskip=0pt
\obeyspaces
i75\ :\ dim\ Proj(R/Lines),\ degree(Lines)\leavevmode\hss\endgraf
\penalty-500\leavevmode\hss\endgraf
o75\ =\ (0,\ 12)\leavevmode\hss\endgraf
\penalty-500\leavevmode\hss\endgraf
o75\ :\ Sequence\leavevmode\hss\endgraf
\endgroup
\penalty-1000
\par
\vskip 1 pt
\noindent

We investigate the rest of the zeroes, which we obtain 
by taking the ideal
quotient of $I$ and the ideal of lines.
As computed above, this has dimension 1 and degree 4.
\par
\vskip 5 pt
\begingroup
\tteight
\baselineskip=10.01pt
\lineskip=0pt
\obeyspaces
i76\ :\ Junk\ =\ I\ :\ Lines;\leavevmode\hss\endgraf
\penalty-500\leavevmode\hss\endgraf
o76\ :\ Ideal\ of\ R\leavevmode\hss\endgraf
\endgroup
\penalty-1000
\par
\vskip 1 pt
\noindent
\par
\vskip 5 pt
\begingroup
\tteight
\baselineskip=10.01pt
\lineskip=0pt
\obeyspaces
i77\ :\ dim\ Proj(R/Junk),\ degree\ Junk\leavevmode\hss\endgraf
\penalty-500\leavevmode\hss\endgraf
o77\ =\ (1,\ 4)\leavevmode\hss\endgraf
\penalty-500\leavevmode\hss\endgraf
o77\ :\ Sequence\leavevmode\hss\endgraf
\endgroup
\penalty-1000
\par
\vskip 1 pt
\noindent
We find the support of this extraneous component by taking its
radical.
\par
\vskip 5 pt
\begingroup
\tteight
\baselineskip=10.01pt
\lineskip=0pt
\obeyspaces
i78\ :\ radical(Junk)\leavevmode\hss\endgraf
\penalty-500\leavevmode\hss\endgraf
\ \ \ \ \ \ \ \ \ \ \ \ \ \ \ \ \ \ \ \ \ \ \ \ \ \ \ \ \ \ \ \ \ \ \ \ \ \ \ \ \ 2\ \ \ \ \ \ \ \ \ 2\ \ \ \ \ \ \ \ \ 2\leavevmode\hss\endgraf
o78\ =\ ideal\ (p\ \ \ \ \ \ ,\ p\ \ \ \ \ \ ,\ p\ \ \ \ \ \ ,\ p\ \ \ \ \ \ \ +\ p\ \ \ \ \ \ \ +\ p\ \ \ \ \ \ )\leavevmode\hss\endgraf
\ \ \ \ \ \ \ \ \ \ \ \ \ \ {\char`\{}0,\ 3{\char`\}}\ \ \ {\char`\{}0,\ 2{\char`\}}\ \ \ {\char`\{}0,\ 1{\char`\}}\ \ \ {\char`\{}1,\ 2{\char`\}}\ \ \ \ {\char`\{}1,\ 3{\char`\}}\ \ \ \ {\char`\{}2,\ 3{\char`\}}\leavevmode\hss\endgraf
\penalty-500\leavevmode\hss\endgraf
o78\ :\ Ideal\ of\ R\leavevmode\hss\endgraf
\endgroup
\penalty-1000
\par
\vskip 1 pt
\noindent
From this, we see that the extraneous component is supported on an imaginary
conic in the ${\mathbb P}^2$ of lines at infinity.
\smallskip

To understand the geometry behind this computation, observe that
the sphere with radius $r$ and center $(a,b,c)$ has homogeneous equation
$$
  (x-wa)^2+(y-wb)^2+(z-wc)^2\ =\ r^2w^2\,.
$$
At infinity, $w=0$, this has equation
$$
  x^2+y^2+z^2\ =\ 0\,.
$$
The extraneous component is supported on the set of tangent
lines to this imaginary conic.
Aluffi\index{Aluffi, P.} and Fulton\index{Fulton, W.}~\cite{SO:AF} studied
this problem, using geometry to 
identify the extraneous ideal and the excess intersection
formula~\cite{SO:FM76} to obtain the answer of 12. 
Their techniques show that there will be 12 isolated lines tangent to 4
quadrics which have a smooth conic in common.

When the quadrics are spheres, the conic is the imaginary conic at infinity.
Fulton\index{Fulton, W.} asked the following question:
Can all 12 lines be real if the (real) four quadrics share a real conic?
We answer his question in the affirmative in the next section.

\subsection{Lines tangent to real quadrics sharing a real conic}
We consider four quadrics in ${\mathbb P}^3_{\mathbb R}$ sharing a
non-singular conic, which we will take to be at infinity so that we may use
local coordinates for ${\bf G}_{2,4}$ in our computations.
The variety ${\mathcal V}(q)\subset{\mathbb P}^3_{\mathbb R}$ of a
nondegenerate quadratic form $q$ is determined up to isomorphism by the
absolute value of the signature\index{bilinear form!signature} $\sigma$ of the
associated bilinear form. 
Thus there are three possibilities, 0, 2, or 4, for $|\sigma|$.

When $|\sigma|=4$, the real quadric ${\mathcal V}(q)$ is empty.
The associated symmetric matrix $M$ is conjugate to the identity
matrix, so $\wedge^2M$ is also conjugate to the identity matrix.
Hence ${\mathcal V}(\wedge^2q)$ contains no real points.
Thus we need not consider quadrics with $|\sigma|=4$.

When $|\sigma|=2$, we have ${\mathcal V}(q)\simeq S^2$, the 2-sphere.
If the conic at infinity is imaginary, then 
${\mathcal V}(q)\subset{\mathbb R}^3$ is an ellipsoid\index{ellipsoid}.
If the conic at infinity is real, then ${\mathcal V}(q)\subset{\mathbb R}^3$ is 
a hyperboloid\index{hyperboloid} of two sheets.
When $\sigma=0$, we have ${\mathcal V}(q)\simeq S^1\times S^1$, a torus.
In this case, ${\mathcal V}(q)\subset{\mathbb R}^3$ is 
a hyperboloid of one sheet and the conic at infinity is real.

Thus either we have 4 ellipsoids sharing an imaginary conic at infinity,
which we studied in Section~\ref{sec:12lines}; or else we have four
hyperboloids sharing a real conic at infinity, and there are five 
possible combinations of hyperboloids of one or two sheets sharing a real
conic at infinity.
This gives six topologically distinct possibilities in all.

\begin{theorem}
For each of the six topologically distinct possibilities of four real
quadrics sharing a smooth conic at infinity, there exist four quadrics
having the property that each of the 12 lines in ${\mathbb C}^3$
simultaneously tangent to the four quadrics is real.
\end{theorem}

\begin{proof}
By the computation in Section~\ref{sec:12lines}, 
we need only check the five possibilities for hyperboloids\index{hyperboloid}.
We fix the conic at infinity to be $x^2+y^2-z^2=0$.
Then the general hyperboloid of two sheets containing this conic has
equation in ${\mathbb R}^3$
\begin{equation}\label{eq:twoSheet}
  (x-a)^2+(y-b)^2-(z-c)^2+r\ =\ 0\,,
\end{equation}
(with $r>0$).
The command {\tt Two(a,b,c,r)} generates the associated 
symmetric matrix.
\par
\vskip 5 pt
\begingroup
\tteight
\baselineskip=10.01pt
\lineskip=0pt
\obeyspaces
i79\ :\ Two\ =\ (a,\ b,\ c,\ r)\ ->\ (\leavevmode\hss\endgraf
\ \ \ \ \ \ \ \ \ \ \ matrix{\char`\{}{\char`\{}a{\char`\^}2\ +\ b{\char`\^}2\ -\ c{\char`\^}2\ +\ r\ ,-a\ ,-b\ ,\ c\ {\char`\}},\leavevmode\hss\endgraf
\ \ \ \ \ \ \ \ \ \ \ \ \ \ \ \ \ \ {\char`\{}\ \ \ \ \ \ \ \ \ -a\ \ \ \ \ \ \ \ \ ,\ 1\ ,\ 0\ ,\ 0\ {\char`\}},\leavevmode\hss\endgraf
\ \ \ \ \ \ \ \ \ \ \ \ \ \ \ \ \ \ {\char`\{}\ \ \ \ \ \ \ \ \ -b\ \ \ \ \ \ \ \ \ ,\ 0\ ,\ 1\ ,\ 0\ {\char`\}},\leavevmode\hss\endgraf
\ \ \ \ \ \ \ \ \ \ \ \ \ \ \ \ \ \ {\char`\{}\ \ \ \ \ \ \ \ \ \ c\ \ \ \ \ \ \ \ \ ,\ 0\ ,\ 0\ ,-1\ {\char`\}}{\char`\}}\leavevmode\hss\endgraf
\ \ \ \ \ \ \ \ \ \ \ );\leavevmode\hss\endgraf
\endgroup
\penalty-1000
\par
\vskip 1 pt
\noindent
The general hyperboloid of one sheet containing the conic
$x^2+y^2-z^2=0$ at infinity has equation in ${\mathbb R}^3$
\begin{equation}\label{eq:oneSheet}
  (x-a)^2+(y-b)^2-(z-c)^2-r\ =\ 0\,,
\end{equation}
(with $r>0$).
The command {\tt One(a,b,c,r)} generates the associated 
symmetric matrix.
\par
\vskip 5 pt
\begingroup
\tteight
\baselineskip=10.01pt
\lineskip=0pt
\obeyspaces
i80\ :\ One\ =\ (a,\ b,\ c,\ r)\ ->\ (\leavevmode\hss\endgraf
\ \ \ \ \ \ \ \ \ \ \ matrix{\char`\{}{\char`\{}a{\char`\^}2\ +\ b{\char`\^}2\ -\ c{\char`\^}2\ -\ r\ ,-a\ ,-b\ ,\ c\ {\char`\}},\leavevmode\hss\endgraf
\ \ \ \ \ \ \ \ \ \ \ \ \ \ \ \ \ \ {\char`\{}\ \ \ \ \ \ \ \ \ -a\ \ \ \ \ \ \ \ \ ,\ 1\ ,\ 0\ ,\ 0\ {\char`\}},\leavevmode\hss\endgraf
\ \ \ \ \ \ \ \ \ \ \ \ \ \ \ \ \ \ {\char`\{}\ \ \ \ \ \ \ \ \ -b\ \ \ \ \ \ \ \ \ ,\ 0\ ,\ 1\ ,\ 0\ {\char`\}},\leavevmode\hss\endgraf
\ \ \ \ \ \ \ \ \ \ \ \ \ \ \ \ \ \ {\char`\{}\ \ \ \ \ \ \ \ \ \ c\ \ \ \ \ \ \ \ \ ,\ 0\ ,\ 0\ ,-1\ {\char`\}}{\char`\}}\leavevmode\hss\endgraf
\ \ \ \ \ \ \ \ \ \ \ );\leavevmode\hss\endgraf
\endgroup
\penalty-1000
\par
\vskip 1 pt
\noindent

We consider $i$ quadrics of two sheets~(\ref{eq:twoSheet}) and $4-i$
quadrics of one sheet~(\ref{eq:oneSheet}).
For each of these cases, the table below displays four 4-tuples of data
$(a,b,c,r)$ which give 12 real lines.
(The data for the hyperboloids of one sheet are listed first.)
\begin{center}
\begin{tabular}{|c|l|}\hline
$\;i$\;&\hspace{13em}Data\\\hline
$0$&$\hspace{.915em}(5,3,3,16),\hspace{2.185em}(5,-4,2,1),\hspace{1.27em}
                             (-3,-1,1,1),\hspace{.88em}(2,-7,0,1)$\\\hline
$1$&$\hspace{.385em}(3,-2,-3,6),\hspace{.885em}(-3,-7,-6,7),\hspace{.5em}
                            (-6,3,-5,2),\hspace{.865em}(1,6,-2,5)$\\\hline
$2$&$\hspace{1.165em}(6,4,6,4),\hspace{2.43em}(-1,3,3,6),\hspace{1.265em}
                             (-7,-2,3,3),\hspace{.5em}(-6,7,-2,5)$\\\hline
$3$&$(-1,-4,-1,1),\hspace{.885em}(-3,3,-1,1),\hspace{.885em}
              \hspace{.4em}(-7,6,2,9),\hspace{1.025em}(5,6,-1,12)$\\\hline
$4$&$\hspace{.525em}(5,2,-1,25),\hspace{1.555em}(6,-6,2,25),\hspace{1.03em}
                 \hspace{.4em}(-7,1,6,1),\hspace{1.65em}(3,1,0,1)$\\\hline
\end{tabular}
\end{center}



We test each of these, using the formulation in local coordinates of
Section~\ref{sec:12lines}. 
\par
\vskip 5 pt
\begingroup
\tteight
\baselineskip=10.01pt
\lineskip=0pt
\obeyspaces
i81\ :\ R\ =\ QQ[y11,\ y12,\ y21,\ y22];\leavevmode\hss\endgraf
\endgroup
\penalty-1000
\par
\vskip 1 pt
\noindent
\par
\vskip 5 pt
\begingroup
\tteight
\baselineskip=10.01pt
\lineskip=0pt
\obeyspaces
i82\ :\ I\ =\ ideal\ (tangentTo(One(\ 5,\ 3,\ 3,16)),\ \leavevmode\hss\endgraf
\ \ \ \ \ \ \ \ \ \ \ \ \ \ \ \ \ \ tangentTo(One(\ 5,-4,\ 2,\ 1)),\ \ \leavevmode\hss\endgraf
\ \ \ \ \ \ \ \ \ \ \ \ \ \ \ \ \ \ tangentTo(One(-3,-1,\ 1,\ 1)),\ \leavevmode\hss\endgraf
\ \ \ \ \ \ \ \ \ \ \ \ \ \ \ \ \ \ tangentTo(One(\ 2,-7,\ 0,\ 1)));\leavevmode\hss\endgraf
\penalty-500\leavevmode\hss\endgraf
o82\ :\ Ideal\ of\ R\leavevmode\hss\endgraf
\endgroup
\penalty-1000
\par
\vskip 1 pt
\noindent
\par
\vskip 5 pt
\begingroup
\tteight
\baselineskip=10.01pt
\lineskip=0pt
\obeyspaces
i83\ :\ numRealSturm(charPoly(promote(y22,\ R/I),\ Z))\leavevmode\hss\endgraf
\penalty-500\leavevmode\hss\endgraf
o83\ =\ 12\leavevmode\hss\endgraf
\endgroup
\penalty-1000
\par
\vskip 1 pt
\noindent
\par
\vskip 5 pt
\begingroup
\tteight
\baselineskip=10.01pt
\lineskip=0pt
\obeyspaces
i84\ :\ I\ =\ ideal\ (tangentTo(One(\ 3,-2,-3,\ 6)),\ \leavevmode\hss\endgraf
\ \ \ \ \ \ \ \ \ \ \ \ \ \ \ \ \ \ tangentTo(One(-3,-7,-6,\ 7)),\ \ \leavevmode\hss\endgraf
\ \ \ \ \ \ \ \ \ \ \ \ \ \ \ \ \ \ tangentTo(One(-6,\ 3,-5,\ 2)),\ \leavevmode\hss\endgraf
\ \ \ \ \ \ \ \ \ \ \ \ \ \ \ \ \ \ tangentTo(Two(\ 1,\ 6,-2,\ 5)));\leavevmode\hss\endgraf
\penalty-500\leavevmode\hss\endgraf
o84\ :\ Ideal\ of\ R\leavevmode\hss\endgraf
\endgroup
\penalty-1000
\par
\vskip 1 pt
\noindent
\par
\vskip 5 pt
\begingroup
\tteight
\baselineskip=10.01pt
\lineskip=0pt
\obeyspaces
i85\ :\ numRealSturm(charPoly(promote(y22,\ R/I),\ Z))\leavevmode\hss\endgraf
\penalty-500\leavevmode\hss\endgraf
o85\ =\ 12\leavevmode\hss\endgraf
\endgroup
\penalty-1000
\par
\vskip 1 pt
\noindent
\par
\vskip 5 pt
\begingroup
\tteight
\baselineskip=10.01pt
\lineskip=0pt
\obeyspaces
i86\ :\ I\ =\ ideal\ (tangentTo(One(\ 6,\ 4,\ 6,\ 4)),\ \ \leavevmode\hss\endgraf
\ \ \ \ \ \ \ \ \ \ \ \ \ \ \ \ \ \ tangentTo(One(-1,\ 3,\ 3,\ 6)),\ \leavevmode\hss\endgraf
\ \ \ \ \ \ \ \ \ \ \ \ \ \ \ \ \ \ tangentTo(Two(-7,-2,\ 3,\ 3)),\ \leavevmode\hss\endgraf
\ \ \ \ \ \ \ \ \ \ \ \ \ \ \ \ \ \ tangentTo(Two(-6,\ 7,-2,\ 5)));\leavevmode\hss\endgraf
\penalty-500\leavevmode\hss\endgraf
o86\ :\ Ideal\ of\ R\leavevmode\hss\endgraf
\endgroup
\penalty-1000
\par
\vskip 1 pt
\noindent
\par
\vskip 5 pt
\begingroup
\tteight
\baselineskip=10.01pt
\lineskip=0pt
\obeyspaces
i87\ :\ numRealSturm(charPoly(promote(y22,\ R/I),\ Z))\leavevmode\hss\endgraf
\penalty-500\leavevmode\hss\endgraf
o87\ =\ 12\leavevmode\hss\endgraf
\endgroup
\penalty-1000
\par
\vskip 1 pt
\noindent

\par
\vskip 5 pt
\begingroup
\tteight
\baselineskip=10.01pt
\lineskip=0pt
\obeyspaces
i88\ :\ I\ =\ ideal\ (tangentTo(One(-1,-4,-1,\ 1)),\leavevmode\hss\endgraf
\ \ \ \ \ \ \ \ \ \ \ \ \ \ \ \ \ \ tangentTo(Two(-3,\ 3,-1,\ 1)),\ \ \leavevmode\hss\endgraf
\ \ \ \ \ \ \ \ \ \ \ \ \ \ \ \ \ \ tangentTo(Two(-7,\ 6,\ 2,\ 9)),\ \leavevmode\hss\endgraf
\ \ \ \ \ \ \ \ \ \ \ \ \ \ \ \ \ \ tangentTo(Two(\ 5,\ 6,-1,12)));\leavevmode\hss\endgraf
\penalty-500\leavevmode\hss\endgraf
o88\ :\ Ideal\ of\ R\leavevmode\hss\endgraf
\endgroup
\penalty-1000
\par
\vskip 1 pt
\noindent
\par
\vskip 5 pt
\begingroup
\tteight
\baselineskip=10.01pt
\lineskip=0pt
\obeyspaces
i89\ :\ numRealSturm(charPoly(promote(y22,\ R/I),\ Z))\leavevmode\hss\endgraf
\penalty-500\leavevmode\hss\endgraf
o89\ =\ 12\leavevmode\hss\endgraf
\endgroup
\penalty-1000
\par
\vskip 1 pt
\noindent
\par
\vskip 5 pt
\begingroup
\tteight
\baselineskip=10.01pt
\lineskip=0pt
\obeyspaces
i90\ :\ I\ =\ ideal\ (tangentTo(Two(\ 5,\ 2,-1,25)),\ \leavevmode\hss\endgraf
\ \ \ \ \ \ \ \ \ \ \ \ \ \ \ \ \ \ tangentTo(Two(\ 6,-6,\ 2,25)),\ \leavevmode\hss\endgraf
\ \ \ \ \ \ \ \ \ \ \ \ \ \ \ \ \ \ tangentTo(Two(-7,\ 1,\ 6,\ 1)),\ \leavevmode\hss\endgraf
\ \ \ \ \ \ \ \ \ \ \ \ \ \ \ \ \ \ tangentTo(Two(\ 3,\ 1,\ 0,\ 1)));\leavevmode\hss\endgraf
\penalty-500\leavevmode\hss\endgraf
o90\ :\ Ideal\ of\ R\leavevmode\hss\endgraf
\endgroup
\penalty-1000
\par
\vskip 1 pt
\noindent
\par
\vskip 5 pt
\begingroup
\tteight
\baselineskip=10.01pt
\lineskip=0pt
\obeyspaces
i91\ :\ numRealSturm(charPoly(promote(y22,\ R/I),\ Z))\leavevmode\hss\endgraf
\penalty-500\leavevmode\hss\endgraf
o91\ =\ 12\leavevmode\hss\endgraf
\endgroup
\penalty-1000
\par
\vskip 1 pt
\noindent
%
\end{proof}

In each of these enumerative problems, we have checked that every possible
number of real solutions (0, 2, 4, 6, 8, 10, or 12) can occur.

\subsection{Generalization to higher dimensions}
We consider lines tangent to quadrics in higher dimensions.
First, we reinterpret the action of $\wedge^rV^*$ on $\wedge^rV$
described in~(\ref{eq:wedge}) as follows.
The vectors ${\bf x}_1,\ldots,{\bf x}_r$ and ${\bf y}\!_1,\ldots,{\bf y}\!_r$ 
define maps $\alpha\colon k^r\to V^*$ and $\beta\colon k^r\to V$.
The matrix $[({\bf x}_i,\,{\bf y}\!_j)]$ is the matrix of the bilinear form
on $k^r$ given by 
$\langle{\bf u},\,{\bf v}\rangle:= (\alpha({\bf u}),\,\beta({\bf v}))$.
Thus~(\ref{eq:wedge}) vanishes when the bilinear form 
$\langle\,\cdot\,,\,\cdot\,\rangle$ on $k^r$ is degenerate.

Now suppose that we have a quadratic form $q$ on $V$ given by a symmetric map 
$\varphi\colon V\to V^*$.
This induces a quadratic form and hence a quadric on any $r$-plane $H$ in $V$ 
(with $H\not\subset{\mathcal V}(q)$).
This induced quadric is singular when $H$ is tangent to ${\mathcal V}(q)$.
Since a quadratic form is degenerate only when the associated projective
quadric is singular, we see that
$H$ is tangent to the quadric
${\mathcal V}(q)$ if and only if 
$(\wedge^r\varphi(\wedge^rH),\,\wedge^rH)=0$.
(This includes the case $H\subset{\mathcal V}(q)$.)
We summarize this argument.

\begin{theorem}
Let $\varphi\colon V\to V^*$ be a linear map with resulting bilinear form
$(\varphi({\bf u}),\,{\bf v})$.
Then the locus of $r$-planes in $V$ for which the restriction of this form
is degenerate is the set of $r$-planes $H$ whose Pl\"ucker
coordinates\index{Pl\"ucker coordinate} are
isotropic, $(\wedge^r\varphi(\wedge^rH),\,\wedge^rH)=0$, with respect to the
induced form on $\wedge^rV$.

When $\varphi$ is symmetric, this is the locus of $r$-planes tangent to the
associated quadric in ${\mathbb P}(V)$.
\end{theorem}

We explore the problem of lines tangent to quadrics in ${\mathbb P}^n$.
From the calculations of Section~\ref{sec:global}, we do not expect this to
be interesting if the quadrics are general.
(This is borne out for ${\mathbb P}^4$:
 we find 320 lines in ${\mathbb P}^4$ tangent to 6 general quadrics.
 This is the B\'ezout number\index{B\'ezout number}, as $\deg{\bf G}_{2,5}=5$
 and the condition to be tangent to a quadric has degree 2.)
This problem is interesting if the quadrics 
in ${\mathbb P}^n$ share a quadric in a ${\mathbb P}^{n-1}$.
We propose studying such enumerative problems, both determining the 
number of solutions for general such quadrics, and investigating whether or
not it is possible to have all solutions be real. 

We use \Mtwo{}\/ to compute the expected number of
solutions to this problem when $r=2$ and $n=4$.
We first define some functions for this computation, which will involve
counting the degree of the ideal of lines in ${\mathbb P}^4$ tangent to 6
general spheres\index{sphere}.
Here, $X$ gives local coordinates for the Grassmannian\index{Grassmannian},
$M$ is a symmetric matrix, {\tt tanQuad} gives the equation in $X$ for the
lines tangent to the quadric given by $M$.
\par
\vskip 5 pt
\begingroup
\tteight
\baselineskip=10.01pt
\lineskip=0pt
\obeyspaces
i92\ :\ tanQuad\ =\ (M,\ X)\ ->\ (\leavevmode\hss\endgraf
\ \ \ \ \ \ \ \ \ \ \ u\ :=\ X{\char`\^}{\char`\{}0{\char`\}};\leavevmode\hss\endgraf
\ \ \ \ \ \ \ \ \ \ \ v\ :=\ X{\char`\^}{\char`\{}1{\char`\}};\leavevmode\hss\endgraf
\ \ \ \ \ \ \ \ \ \ \ (u\ *\ M\ *\ transpose\ v){\char`\^}2\ -\ \leavevmode\hss\endgraf
\ \ \ \ \ \ \ \ \ \ \ (u\ *\ M\ *\ transpose\ u)\ *\ (v\ *\ M\ *\ transpose\ v)\leavevmode\hss\endgraf
\ \ \ \ \ \ \ \ \ \ \ );\leavevmode\hss\endgraf
\endgroup
\penalty-1000
\par
\vskip 1 pt
\noindent
{\tt nSphere} gives the matrix $M$ for a sphere with
center {\tt V} and squared radius {\tt r}, and {\tt V} and {\tt r} give random
data for a sphere.
\par
\vskip 5 pt
\begingroup
\tteight
\baselineskip=10.01pt
\lineskip=0pt
\obeyspaces
i93\ :\ nSphere\ =\ (V,\ r)\ ->\ \leavevmode\hss\endgraf
\ \ \ \ \ \ \ \ \ \ \ \ \ \ \ (matrix\ {\char`\{}{\char`\{}r\ +\ V\ *\ transpose\ V{\char`\}}{\char`\}}\ ||\ transpose\ V\ )\ |\leavevmode\hss\endgraf
\ \ \ \ \ \ \ \ \ \ \ \ \ \ \ (\ V\ ||\ id{\char`\_}((ring\ r){\char`\^}n)\leavevmode\hss\endgraf
\ \ \ \ \ \ \ \ \ \ \ \ \ \ \ );\leavevmode\hss\endgraf
\endgroup
\penalty-1000
\par
\vskip 1 pt
\noindent
\par
\vskip 5 pt
\begingroup
\tteight
\baselineskip=10.01pt
\lineskip=0pt
\obeyspaces
i94\ :\ V\ =\ ()\ ->\ matrix\ table(1,\ n,\ (i,j)\ ->\ random(0,\ R));\leavevmode\hss\endgraf
\endgroup
\penalty-1000
\par
\vskip 1 pt
\noindent
\par
\vskip 5 pt
\begingroup
\tteight
\baselineskip=10.01pt
\lineskip=0pt
\obeyspaces
i95\ :\ r\ =\ ()\ ->\ random(0,\ R);\leavevmode\hss\endgraf
\endgroup
\penalty-1000
\par
\vskip 1 pt
\noindent
We construct the ambient ring, local coordinates, and the ideal of the
enumerative problem of lines in ${\mathbb P}^4$ tangent to 6 random spheres.
\par
\vskip 5 pt
\begingroup
\tteight
\baselineskip=10.01pt
\lineskip=0pt
\obeyspaces
i96\ :\ n\ =\ 4;\leavevmode\hss\endgraf
\endgroup
\penalty-1000
\par
\vskip 1 pt
\noindent
\par
\vskip 5 pt
\begingroup
\tteight
\baselineskip=10.01pt
\lineskip=0pt
\obeyspaces
i97\ :\ R\ =\ ZZ/1009[flatten(table(2,\ n-1,\ (i,j)\ ->\ z{\char`\_}(i,j)))];\leavevmode\hss\endgraf
\endgroup
\penalty-1000
\par
\vskip 1 pt
\noindent
\par
\vskip 5 pt
\begingroup
\tteight
\baselineskip=10.01pt
\lineskip=0pt
\obeyspaces
i98\ :\ X\ =\ 1\ |\ matrix\ table(2,\ n-1,\ (i,j)\ ->\ z{\char`\_}(i,j))\leavevmode\hss\endgraf
\penalty-500\leavevmode\hss\endgraf
o98\ =\ |\ 1\ 0\ z{\char`\_}(0,0)\ z{\char`\_}(0,1)\ z{\char`\_}(0,2)\ |\leavevmode\hss\endgraf
\ \ \ \ \ \ |\ 0\ 1\ z{\char`\_}(1,0)\ z{\char`\_}(1,1)\ z{\char`\_}(1,2)\ |\leavevmode\hss\endgraf
\penalty-500\leavevmode\hss\endgraf
\ \ \ \ \ \ \ \ \ \ \ \ \ \ 2\ \ \ \ \ \ \ 5\leavevmode\hss\endgraf
o98\ :\ Matrix\ R\ \ <---\ R\leavevmode\hss\endgraf
\endgroup
\penalty-1000
\par
\vskip 1 pt
\noindent
\par
\vskip 5 pt
\begingroup
\tteight
\baselineskip=10.01pt
\lineskip=0pt
\obeyspaces
i99\ :\ I\ =\ ideal\ (apply(1..(2*n-2),\ \leavevmode\hss\endgraf
\ \ \ \ \ \ \ \ \ \ \ \ \ \ \ \ \ \ \ \ \ i\ ->\ tanQuad(nSphere(V(),\ r()),\ X)));\leavevmode\hss\endgraf
\penalty-500\leavevmode\hss\endgraf
o99\ :\ Ideal\ of\ R\leavevmode\hss\endgraf
\endgroup
\penalty-1000
\par
\vskip 1 pt
\noindent
We find there are 24 lines in ${\mathbb P}^4$ tangent to 6 general spheres.
\par
\vskip 5 pt
\begingroup
\tteight
\baselineskip=10.01pt
\lineskip=0pt
\obeyspaces
i100\ :\ dim\ I,\ degree\ I\leavevmode\hss\endgraf
\penalty-500\leavevmode\hss\endgraf
o100\ =\ (0,\ 24)\leavevmode\hss\endgraf
\penalty-500\leavevmode\hss\endgraf
o100\ :\ Sequence\leavevmode\hss\endgraf
\endgroup
\penalty-1000
\par
\vskip 1 pt
\noindent
The expected numbers of solutions we have obtained in this way are displayed
in the table below.
The numbers in boldface are those which are proven.
\begin{center}
\begin{tabular}{|c|c|c|c|c|c|}\hline
  $n$&2&3&4&5&6\\\hline
  \# expected\;&\;{\bf 4}\;&\;{\bf 12}\;&\;24\;&\;48\;&\;96\;\\\hline
\end{tabular}
\end{center}

\section*{Acknowledgements}
We thank Dan Grayson\index{Grayson, D.} and 
Bernd Sturmfels\index{Sturmfels, B.};
Some of the procedures in this chapter were written by Dan Grayson
and the calculation in Section 5.2 is due to Bernd Sturmfels.

\begin{theindex}

  \item \Mtwo{}, 3, 4, 7, 15, 18

  \indexspace

  \item Aluffi, P., 23
  \item artinian, {\it see also} ideal, zero-dimensional, 2

  \indexspace

  \item B\'ezout number, 1, 21, 26
  \item B\'ezout Theorem, 1, 11, 17
  \item bilinear form, 21
    \subitem signature, 8, 23
    \subitem symmetric, 8
  \item Bottema, O., 12

  \indexspace

  \item cylinder, 12

  \indexspace

  \item discriminant, 13

  \indexspace

  \item eliminant, 4--6, 14
  \item elimination theory, 4
  \item ellipsoid, 23
  \item enumerative geometry, 1, 3, 4, 10, 11, 14
    \subitem real, 7, 18
  \item enumerative problem, 3, 11, 14, 20
    \subitem fully real, 17, 19
  \item Eremenko, A., 18

  \indexspace

  \item field
    \subitem algebraically closed, 1, 13
    \subitem splitting, 4
  \item Fulton, W., 23

  \indexspace

  \item Gabrielov, A., 18
  \item Gr\"obner basis, 4
    \subitem reduced, 15
  \item Grassmannian, 12--14, 26
    \subitem local coordinates, 16
    \subitem not a complete intersection, 16
  \item Grayson, D., 27

  \indexspace

  \item Hilbert function, 15
  \item homotopy
    \subitem B\'ezout, 11
    \subitem Gr\"obner, 11, 17
    \subitem optimal, 11
  \item homotopy continuation, 10
  \item hyperboloid, 23, 24

  \indexspace

  \item ideal
    \subitem degree, 2, 3
    \subitem radical, 2, 4, 19
    \subitem zero-dimensional, 2, 18, 22
  \item initial ideal
    \subitem square-free, 11, 17

  \indexspace

  \item Lichtblau, D., 12

  \indexspace

  \item Macdonald, I., 13

  \indexspace

  \item Pach, J., 13
  \item Pl\"ucker coordinate, 15, 20, 21, 26
  \item Pl\"ucker embedding, 15
  \item Pl\"ucker ideal, 15
  \item polynomial equations, 1, 3
    \subitem deficient, 1
    \subitem overdetermined, 16

  \indexspace

  \item quadratic form, 13, 21

  \indexspace

  \item rational normal curve, 17, 18

  \indexspace

  \item saturate, 3, 22
  \item Schubert calculus, 14, 17--19
  \item Shape Lemma, 4, 19
  \item Shapiro, B., 18
  \item Shapiro, M., 18
  \item Shapiros's Conjecture, 18
  \item solving polynomial equations, 3
    \subitem real solutions, 6
    \subitem via eigenvectors, 6
    \subitem via elimination, 4
    \subitem via numerical homotopy, 10
  \item sphere, 13, 22, 26
  \item Stickelberger's Theorem, 6
  \item Sturm sequence, 7
  \item Sturmfels, B., 27
  \item symbolic computation, 2

  \indexspace

  \item Theobald, Th., 13
  \item trace form, 8

  \indexspace

  \item Veldkamp, G., 12

\end{theindex}


\begin{thebibliography}{10}

\bibitem{SO:AF}
P.~Aluffi and W.~Fulton, 2000.
\newblock Private Communication.

\bibitem{SO:BMMT}
E.~Becker, M.~G. Marinari, T.~Mora, and C.~Traverso.
\newblock The shape of the {S}hape {L}emma.
\newblock In {\em Proceedings ISSAC-94}, pages 129--133, 1993.

\bibitem{SO:BW}
E.~Becker and Th. W{\"o}ermann.
\newblock On the trace formula for quadratic forms.
\newblock In {\em Recent advances in real algebraic geometry and quadratic
  forms (Berkeley, CA, 1990/1991; San Francisco, CA, 1991)}, pages 271--291.
  Amer. Math. Soc., Providence, RI, 1994.

\bibitem{SO:Bernstein}
D.~N. Bernstein.
\newblock The number of roots of a system of equations.
\newblock {\em Funct. Anal. Appl.}, 9:183--185, 1975.

\bibitem{SO:BV77}
O.~Bottema and G.R. Veldkamp.
\newblock On the lines in space with equal distances to $n$ given points.
\newblock {\em Geometrie Dedicata}, 6:121--129, 1977.

\bibitem{SO:CCS}
A.~M. Cohen, H.~Cuypers, and H.~Sterk, editors.
\newblock {\em Some Tapas of Computer Algebra}.
\newblock Springer-Varlag, 1999.

\bibitem{SO:CLO92}
D.~Cox, J.~Little, and D.~O'Shea.
\newblock {\em Ideals, Varieties, Algorithms: An Introduction to Computational
  Algebraic Geometry and Commutative Algebra}.
\newblock UTM. Springer-Verlag, New York, 1992.

\bibitem{SO:DeKh00}
A.~Degtyarev and V.~Kharlamov.
\newblock Topological properties of real algebraic curves: du c\^{o}t\'e de
  chez {R}okhlin.
\newblock www.arXiv.org/math.AG/0004134, 2000.

\bibitem{SO:MR97a:13001}
D.~Eisenbud.
\newblock {\em Commutative Algebra With a View Towards Algebraic Geometry}.
\newblock Number 150 in GTM. Springer-Verlag, 1995.

\bibitem{SO:EG00}
A.~Eremenko and A.~Gabrielov.
\newblock Rational functions with real critical points and {B}.~and
  {M}.~{S}hapiro conjecture in real enumerative geometry.
\newblock MSRI preprint 2000-002.

\bibitem{SO:Fu84a}
W.~Fulton.
\newblock {\em Intersection Theory}.
\newblock Number~2 in Ergebnisse der Math. Springer-Verlag, 1984.

\bibitem{SO:FM76}
W.~Fulton and R.~MacPherson.
\newblock Intersecting cycles on an algebraic variety.
\newblock In P.~Holm, editor, {\em Real and Complex Singularities}, pages
  179--197. Oslo, 1976, Sijthoff and Noordhoff, 1977.

\bibitem{SO:HSS}
B.~Huber, F.~Sottile, and B.~Sturmfels.
\newblock Numerical {S}chubert calculus.
\newblock {\em J.~Symb.~Comp.}, 26(6):767--788, 1998.

\bibitem{SO:MR50:13063}
S.~Kleiman.
\newblock The transversality of a general translate.
\newblock {\em Compositio Math.}, 28:287--297, 1974.

\bibitem{SO:MR48:2152}
S.~Kleiman. and D.~Laksov.
\newblock Schubert calculus.
\newblock {\em Amer. Math. Monthly}, 79:1061--1082, 1972.

\bibitem{SO:Li00}
D.~Lichtblau.
\newblock Finding cylinders through 5 points in $\mathbb{R}^3$.
\newblock mss., danl@wolfram.com, 2000.

\bibitem{SO:MPT00}
I.G. Macdonald, J.~Pach, and Th. Theobald.
\newblock Common tangents to four unit balls in $\mathbb{R}^3$.
\newblock Discr.~Comput.~Geom., to appear, 2000.

\bibitem{SO:PRS}
P.~Pedersen, M.-F. Roy, and A.~Szpirglas.
\newblock Counting real zeros in the multivariate case.
\newblock In {\em Computational {A}lgebraic {G}eometry (Nice, 1992)}, pages
  203--224. Birkh\"auser Boston, Boston, MA, 1993.

\bibitem{SO:RS98}
J.~Rosenthal and F.~Sottile.
\newblock Some remarks on real and complex output feedback.
\newblock {\em Systems Control Lett.}, 33(2):73--80, 1998.
\newblock For a description of the computational aspects, see {\tt
  http://www.nd.edu/\~{}rosen/pole/}.

\bibitem{SO:So97a}
F.~Sottile.
\newblock Enumerative geometry for the real {G}rassmannian of lines in
  projective space.
\newblock {\em Duke Math. J.}, 87(1):59--85, 1997.

\bibitem{SO:So99a}
F.~Sottile.
\newblock The special {S}chubert calculus is real.
\newblock {\em ERA of the AMS}, 5:35--39, 1999.

\bibitem{SO:So_shap-www}
F.~Sottile.
\newblock The conjecture of {S}hapiro and {S}hapiro.
\newblock An archive of computations and computer algebra scripts, {\tt
  http://www.expmath.org/extra/9.2/sottile}, 2000.

\bibitem{SO:So00b}
F.~Sottile.
\newblock Real {S}chubert calculus: Polynomial systems and a conjecture of
  {S}hapiro and {S}hapiro.
\newblock {\em Exper.~Math.}, 9:161--182, 2000.

\bibitem{SO:So_flags}
F.~Sottile.
\newblock Some real and unreal enumerative geometry for flag manifolds.
\newblock {\em Mich. Math. J}, 48:573--592, 2000.

\bibitem{SO:So_trans}
F.~Sottile.
\newblock Elementary transversality in the schubert calculus in any
  characteristic.
\newblock math.AG/0010319, 2000.

\bibitem{SO:Sturmfels_GBCP}
Bernd Sturmfels.
\newblock {\em Gr{\"o}bner Bases and Convex Polytopes}, volume~8 of {\em
  University Lecture Series}.
\newblock American Math. Soc., Providence, RI, 1996.

\bibitem{SO:Ver99}
J.~Verschelde.
\newblock Polynomial homotopies for dense, sparse, and determinantal systems.
\newblock MSRI preprint 1999-041, 1999.

\bibitem{SO:Ver00}
J.~Verschelde.
\newblock Numerical evidence of a conjecture in real algebraic geometry.
\newblock {\em Exper.~Math.}, 9:183--196, 2000.

\end{thebibliography}
\end{document}